\documentclass[letterpaper, 10 pt, conference]{ieeeconf}  %

\IEEEoverridecommandlockouts                              %

\newtheorem{definition}{\bf Definition}

 \newtheorem{prop}{\bf Proposition}

\usepackage{amsfonts,amssymb,amsmath}  
\usepackage{mathtools}
\usepackage{bbm}
\usepackage{breqn}
\usepackage{subcaption}
\usepackage{diagbox}
\usepackage{cite}
\usepackage{hyperref}
\usepackage{cleveref}
\crefname{figure}{Fig.}{Figs.} %
\crefrangeformat{equation}{#3(#1)#4--#5(#2)#6} %
\usepackage{graphicx}
\usepackage{caption}
\usepackage{mathtools}

\usepackage{tikz} 
\usepackage{pgfplots}
\usetikzlibrary{matrix}
\usepgfplotslibrary{fillbetween} %
\pgfplotsset{compat=newest}
\usetikzlibrary{plotmarks}
\usetikzlibrary{arrows.meta}
\usepgfplotslibrary{patchplots}
\usetikzlibrary{arrows,decorations.markings}
\usetikzlibrary{decorations.pathreplacing}
\usetikzlibrary{automata,arrows,positioning,calc}

\title{\LARGE \bf
Deceptive Planning for Resource Allocation
}

\author{Shenghui Chen, Yagiz Savas, Mustafa O. Karabag, Brian M. Sadler, Ufuk Topcu \thanks{S. Chen, Y. Savas, M. O. Karabag, and U. Topcu are with the University of Texas at Austin, TX, USA. E-mails: \{shenghui.chen, yagiz.savas, karabag, utopcu\}@utexas.edu }
\thanks{B. M. Sadler is with the U.S. Army Research Laboratory, MD, USA. E-mail: brian.m.sadler6.civ@army.mil }}

\begin{document}

\maketitle

\begin{abstract}
We consider a team of autonomous agents that navigate in an adversarial environment and aim to achieve a task by allocating their resources over a set of target locations. An adversary in the environment observes the autonomous team's behavior to infer their objective and responds against the team. 
In this setting, we propose strategies for controlling the density of the autonomous team so that they can deceive the adversary regarding their objective while achieving the desired final resource allocation. 
We first develop a prediction algorithm based on the principle of maximum entropy to express the team's behavior expected by the adversary. Then, by measuring the deceptiveness via Kullback-Leibler divergence, we devise convex optimization-based planning algorithms that deceive the adversary by either exaggerating the behavior towards a decoy allocation strategy or creating ambiguity regarding the final allocation strategy. 
A user study with $320$ participants demonstrates that the proposed algorithms are effective for deception and reveal the inherent biases of participants towards proximate goals.
\end{abstract}

\section{Introduction}
In many scenarios, a team of autonomous agents needs to accomplish a task in an adversarial environment. 
Consider a swarm of autonomous drones tasked with securing a region and conducting surveillance missions amidst potential intruders~\cite{saska2016swarm}, or autonomous military robots aiming to control strategic locations from opposing parties on battlefields~\cite{lin2008autonomous}.
Operating in such hostile environments often leads to the team inadvertently leaking critical information, enabling the adversary to devise counter-strategies that thwart task completion.

\begin{figure}[t]
\centering
\resizebox{1\linewidth}{!}{
\begin{tikzpicture}

\node[] at (5,5.5) {\Large \textbf{Adversary prediction model}};
\draw[draw=black, fill = lightgray, opacity = 0.4, rounded corners](2.5,0.5) rectangle ++(5,4.5);

\draw[draw=black, fill = orange, opacity = 0.4, rounded corners](2.8,3) rectangle ++(4.4,1.5);
\node[text width = 4.5 cm,align = center] at (5,4.1) {\large Goal-directed behavior };

\node[text width = 3.5 cm,align = center] at (5,3.5) {\large $c$ $\colon$ $\mathcal{S}$$\times$$\mathcal{A}$$\rightarrow$$[0,\infty)$};

\draw[draw=black, fill = orange, opacity = 0.4,rounded corners](2.8,1) rectangle ++(4.4,1.5);
\node[text width = 4 cm,align = center] at (5,2.1) {\large Expected inefficiency };
\node[text width = 4 cm,align = center] at (5,1.5) {\large $\beta$$\in$$(0,\infty)$};

\draw[->, thick] (7.5,2.75) -- (11-2,2.75);

\node[] at (13.5-2,5.5) {\Large \textbf{Team's planning model}};
\draw[draw=black, fill = lightgray, opacity = 0.4,rounded corners](11-2,0.5) rectangle ++(5,4.5);
\draw[draw=black, fill = orange, opacity = 0.4,rounded corners](11.3-2,3) rectangle ++(4.4,1.5);
\node[align = center] at (13.5-2,4.2) {\large True goal distribution};
\node[align = center] at (13.5-2,3.5) {\large $\boldsymbol{\sigma}_{\star}$$\in$$\mathcal{U}$};

\draw[draw=black, fill = orange, opacity = 0.4,rounded corners](11.3-2,1) rectangle ++(4.4,1.5);
\node[text width = 4 cm,align = center] at (13.5-2,2.1) {\large Type of deception};
\node[text width = 3 cm, align = center] at (13.5-2,1.5) {\large $\text{Deceptiveness}(\mathcal{U}, \pi)$};

\draw[->, thick] (16-2,2.75) -- (19-4,2.75);
\draw[draw=black, fill = lightgray, opacity = 0.4,rounded corners](19-4,1.75) rectangle ++(3,2);
\node[text width = 2 cm, align = center] at (20.35-4,3.2) {\large Deceptive policy};
\node[text width = 2 cm, align = center] at (20.15-4,2.3) {\large  $\pi^{\star}$$\colon$$\mathcal{S}$$\times$$\mathcal{A}$$\rightarrow$$[0,1]$};

\end{tikzpicture}}
\caption{The proposed deceptive resource allocation approach. For a goal distribution (allocation), the adversary prediction model describes how the adversary expects the autonomous team to achieve its final distribution. Based on the predicted policies, the team generates a deceptive policy either exaggerating the team's behavior toward a decoy goal distribution or creating ambiguity regarding the true goal distribution while achieving the desired final allocation.}
\label{flowchart}
\vspace{-10pt} %
\end{figure}
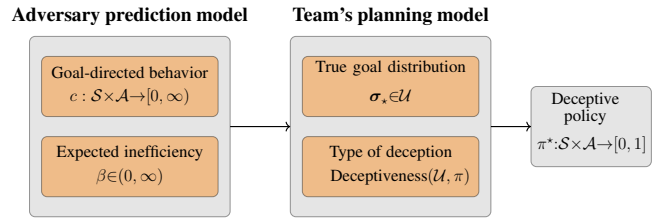

In this paper, we present a systematic approach for a team of agents to complete their task while managing information leakage through deliberate deception.
We consider a setting in which a team consisting of a large number of autonomous agents distribute their resources, i.e., team members, to certain goal locations in an environment. 
Knowing only that the true goal distribution is among a limited set of distributions, an adversary observes the team's behavior in the environment to deduce the team's goal distribution and respond against the team.
In this setting, we develop a swarm control strategy for the autonomous team to allocate their resources to desired locations in a way to deceive the adversary regarding the true distribution. The approach is summarized in Fig. \ref{flowchart}.

We model the prior prediction of the adversary via the maximum entropy principle~\cite{jaynes1957information, ziebart2008maximum}. Specifically, inspired by the experimental studies from the psychology literature~\cite{gergely1995taking}, we assume that the adversary expects the autonomous team to reach their final allocation in the environment through the shortest paths with a certain degree of inefficiency. We generate the expected behavior by solving a constrained optimization problem that combines a cost minimization objective with an entropy regularization.

We model the behavior of the team in the environment as a Markov decision process (MDP). MDPs model sequential decision-making problems under uncertainty and have been widely used to control the high-level behavior of autonomous agents in various applications~\cite{feinberg2012handbook, Puterman, dolgov2006resource,hibbard2020minimizing, witwicki2017autonomous}. We utilize MDPs to synthesize a strategy that controls the density of the team members in the environment while they progress toward achieving the desired final resource allocation. 
Specifically, we synthesize a strategy that maximizes the deceptiveness of the transient behavior while guaranteeing the attainment of the intended final allocation.

We quantify the deceptiveness of the team's behavior as a function of the statistical distance between the observed behavior and the behavior expected to achieve the true objective. In particular, we consider two types of deception, namely, exaggeration and ambiguity, and show how Kullback-Leibler divergence between certain distributions can be used to develop deception metrics. 

This paper has three main contributions. First, we show that an entropy-regularized cost minimization problem in MDPs subject to multiple probabilistic constraints can be formulated as a convex optimization problem and solved efficiently via off-the-shelf solvers. Unlike the existing literature on deception that typically focuses on a single agent with a single reachability objective, this work enables the modeling of adversary predictions in scenarios that involve a swarm of agents with multiple reachability objectives. Second, we introduce novel metrics to quantify the deceptiveness of the team's behavior and present efficient convex optimization-based algorithms to synthesize strategies that control the density of the team and yield globally optimal deceptive behaviors while satisfying multiple reachability constraints. Third, we validate the deceptiveness of the synthesized strategies via a user study with $320$ participants. Our results show that the proposed deceptive algorithms are effective and reveal the inherent biases of participants towards goals closer to the starting point.

\noindent \textbf{Related work:} Several lines of work are related to the deception problem considered in this paper. The most closely related ones are the authors' previous work on supervisory control \cite{karabag2021deception} and deception under uncertainty \cite{savas2021deceptive}. 
The former studies how to deceive a supervisor who provides a reference policy for the agent to follow, inspiring our use of hypothesis testing theory to formulate different deception techniques. The latter considers a deception problem in a single-agent setting, presenting a maximum-entropy-based algorithm for prediction and a linear-programming-based algorithm to optimally reach a single goal. Unlike \cite{savas2021deceptive}, we develop a convex optimization-based prediction algorithm that incorporates reachability objectives for multiple goals. Rather than defining a cost function based on prediction probabilities, we employ Kullback-Leibler divergence to quantify deception.

There is a large body of literature on single-agent deception problems in which the agent aims to reach its goal while deceiving outside observers. The paper \cite{dragan2015deceptive} presents a gradient-descent-based approach to synthesize locally optimal deceptive strategies for reaching a single goal in deterministic environments. They consider both exaggeration and ambiguity types of deception by quantifying deceptiveness as a function of prediction probabilities. 
Unlike \cite{dragan2015deceptive}, we quantify deceptiveness as a function of the statistical distances between density distributions. 
The paper \cite{masters2017deceptive} introduces the notion of the last deceptive point in an environment and presents heuristic approaches to synthesize deceptive policies based on this notion. Similarly, \cite{ornik2018deception} develops deceptive strategies by modeling the observer predictions as a stochastic transition system over potential goals. 
Although the techniques presented in these works are quite insightful, they are designed for single-objective scenarios and are not applicable to situations requiring a team to allocate a certain fraction of their resources over multiple targets.

Game-theoretic approaches are also commonly used to develop deception strategies in various applications. In \cite{li2022dynamic, kulkarni2020deceptive, kulkarni2020decoy}, the authors develop several algorithms to utilize decoys for deception in hypergames. Unlike the efficient algorithms presented in this paper, hypergame formulations, in general, yield computationally intractable solutions that can hardly be applied to large-scale systems. 
Other research \cite{nguyen2019deception} and \cite{wagner2011acting} develop deceptive strategies for specific game types, focusing on finitely repeated and single-stage games, respectively.
Our work differs from these papers as we consider a dynamic system model where an autonomous team needs to navigate in an environment to achieve an objective eventually.

Finally, the literature on goal recognition is also closely related to deception. In \cite{ramirez2010probabilistic,ramirez2011goal,shvo2020active}, the authors develop several algorithms for observers to infer an agent's goal based on its past behavior. These algorithms typically focus on deterministic environments and assume that the agent aims to reach one of the finitely many goals. Since we model the team's behavior as an MDP, the inference techniques presented in this paper also apply to stochastic environments. Moreover, unlike the existing algorithms on goal recognition, the proposed maximum-entropy-based approach can handle scenarios in which the team aims to reach multiple goals with associated probabilities.
\section{Preliminaries}

\noindent\textbf{Notation:} For a given set $\mathcal{S}$, we denote its cardinality by $\lvert \mathcal{S} \rvert$. We define $\mathbb{N}:=\{1,2,3,\ldots\}$, $\mathbb{R}:=(-\infty, \infty)$, and $\mathbb{R}_{\geq 0}:=[0,\infty)$. For a matrix $M\in\mathbb{R}^{n\times m}$, we denote its $(i,j)$-th element by $M_{i,j}$ and its transpose by $M^T$. Finally, for a constant $K\in\mathbb{N}$, we denote the set $\{1,2,\ldots,K\}$ by $[K]$.

\subsection{Markov decision processes}
We consider a team consisting of a large number of autonomous agents. We model the behavior of the team in a stochastic environment with a Markov decision process.

{\setlength{\parindent}{0cm}
\begin{definition}\label{def:MDP}
    A \textit{Markov decision process} (MDP) is a tuple $\mathbb{M}=(\mathcal{S},\alpha,\mathcal{A},P)$ where $\mathcal{S}$ is a finite set of states, $\alpha:\mathcal{S}\rightarrow[0,1]$ is an initial state distribution, $\mathcal{A}$ is a finite set of actions, and $P:\mathcal{S}\times\mathcal{A}\rightarrow[0,1]$ is a transition probability function such that $\sum_{s'\in\mathcal{S}}P(s,a,s')=1$ for all $s\in\mathcal{S}$ and $a\in\mathcal{A}(s)$, where $\mathcal{A}(s)\subseteq\mathcal{A}$ is the set of available actions in state $s$.
\end{definition}}

For notational convenience, we denote the transition probability $P(s,a,s')$ by $P_{s,a,s'}$. A state $s\in\mathcal{S}$ is said to be \textit{absorbing} if $P_{s,a,s}=1$ for all $a\in\mathcal{A}(s)$. 

We control the temporal density of the team through a policy to be applied by every member of the team.

{\setlength{\parindent}{0cm}
\noindent \begin{definition}
For an MDP $\mathbb{M}$, a \textit{policy} $\pi: \mathcal{S}\times\mathcal{A}\rightarrow[0,1]$ is a mapping such that $\sum_{a\in \mathcal{A}(s)}\pi(s,a)=1$ for all $s\in\mathcal{S}$. We denote the set of all policies by $\Pi(\mathbb{M})$. 
\end{definition}}

We note that a policy $\pi$ is traditionally referred to as a \textit{stationary} policy \cite{Puterman}. Although it is possible to consider more general policy classes, the set $\Pi(\mathbb{M})$ is sufficient without loss of generality for the purposes of this paper.

A \textit{path} is a sequence $\varrho=s_1a_1s_2a_2s_3\ldots$ of states and actions which satisfies that $\alpha(s_1)>0$ and $P_{s_t,a_t,s_{t+1}}>0$ for all $t\in\mathbb{N}$. We define the set of all paths in $\mathbb{M}$ with initial distribution $\alpha$ generated under the policy $\pi$ by $Paths^{\pi,\alpha}_{\mathbb{M}}$ and use the standard probability measure over the set $Paths^{\pi,\alpha}_{\mathbb{M}}$~\cite{Model_checking}. Let $\varrho[t]:=s_t$ denote the state visited at the $t$-th step along $\varrho$. For a given state $s\in\mathcal{S}$, we define
\begin{align*}
     \text{Pr}^{\pi}(\alpha \models \lozenge s):=\text{Pr}\left\{\varrho\in Paths^{\pi,\alpha}_{\mathbb{M}}: \exists t\in \mathbb{N}, \varrho[t] = s \right\}
\end{align*} 
as the probability with which the paths generated in $\mathbb{M}$ with initial distribution $\alpha$ under $\pi$ reaches the state $s\in \mathcal{S}$.

\section{Problem Statement}
We consider a team of autonomous agents that are distributed in a stochastic environment. The team aims to navigate through the environment and achieve a desired final distribution, expressing the optimal allocation of resources to certain goal locations. There is an adversary that observes the team's behavior and aims to predict their final distribution to respond against the team's allocation. We note that the distribution of the resources could be an outcome of an underlying game (e.g., a Colonel Blotto game~\cite{roberson2006colonel}, a zero-sum matrix game~\cite{raghavan1994zero}) between the team and the adversary. We study the problem of generating a swarm-control policy for the team members so that they deceive the adversary regarding their final distribution for as long as possible while eventually achieving the desired final distribution. 

Formally, let $\mathcal{U}=\{\boldsymbol{\sigma}_{1},\ldots,\boldsymbol{\sigma}_{N}\}$ be a finite set of goal distributions where $\boldsymbol{\sigma}_{*}$ is the \textit{true goal distribution}. If the adversary knew the true goal distribution, then deception would not be possible as the adversary would respond optimally regardless of the team's behavior in the environment. However, deception becomes possible when the adversary only knows that the true goal distribution belongs to the set $\mathcal{U}$ of potential goal distributions. 

For a given set $\mathcal{U}$ of potential goal distribution and a policy $\pi\in\Pi(\mathbb{M})$, let $\textit{Deceptiveness}(\mathcal{U},\pi)$ be a measure that quantifies the deceptiveness of the team's behavior. In this paper, we aim to synthesize a policy $\pi^{\star}$ such that
\begin{subequations}
\begin{flalign}\label{opt_main_1}
   \pi^{\star} \in \arg \max_{\pi\in\Pi(\mathbb{M})} &\ \ \textit{Deceptiveness}(\mathcal{U},\pi)&& \\
    \text{subject to:} & \ \ \text{Pr}^{\pi}(\alpha \models \lozenge g_i) = \sigma_{*}(g_i) \  \text{for all} \ g_i\in\mathcal{G}.&&\label{opt_main_2}
\end{flalign}
\end{subequations}
The above problem aims to enable the team to deceive the adversary regarding their true objective while guaranteeing that they achieve the final distribution $\boldsymbol{\sigma}_{\star}$. In what follows, we discuss how to formally define the deceptiveness of a team's induced behavior and develop algorithms to solve the problem in \crefrange{opt_main_1}{opt_main_2}. 

Throughout the paper, we assume that the problem in \crefrange{opt_main_1}{opt_main_2} is feasible. The validity of this assumption for a given problem instance can be verified efficiently by solving a linear program, as described in \cite{Marta}. We note that the feasibility assumption holds in many practical settings. For example, it holds when the MDP model has deterministic transitions and there is at least one path from the initial state to each goal state. Some problem instances violate the assumption due to an unachievable goal state distribution. We do not consider these instances as deception is the main focus of this paper.

\section{Expressing Predictions Through The Principle of Maximum Entropy} \label{sec:maxentropy}
To deceive the adversary about the goal distribution, the team needs to know how they associate the observed behavior with a goal distribution. In this section, we introduce an inference model based on the principle of maximum entropy, which characterizes the adversary's predictions.

Experimental studies show that observers typically expect an agent's behavior to be goal-directed and efficient \cite{gergely1995taking}. This expectation can be expressed through the principle of maximum entropy, which prescribes a probability distribution that is “maximally noncommittal with regard to missing information” \cite{jaynes1957information}. In particular, suppose that the adversary believes that the true goal distribution of the team is $\boldsymbol{\sigma}_{i}$. Then, the team's expected behavior for achieving the final distribution $\boldsymbol{\sigma}_{i}$ is described by a policy $\overline{\pi}_i\in\Pi(\mathbb{M})$ such that
\begin{subequations} \label{entropy_regularized_problem}
\begin{align}
    \overline{\pi}_i \in \arg\min_{\pi\in \Pi(\mathbb{M})} &\; \mathbb{E}^{\pi}\left[\sum_{t=1}^{\infty} \Big(c(s_t,a_t) - \beta H(\pi(\cdot | s_t))\Big)\right] \label{entropy_regularized_obj} \\
    \text{subject to: } &\; \text{Pr}^{\pi}(\alpha \models \lozenge g_i) = \sigma_{i}(g_i) \ \text{for all } g_i\in\mathcal{G}. \label{entropy_regularized_cons}
\end{align}
\end{subequations}
In \eqref{entropy_regularized_problem}, $c:\mathcal{S}\times\mathcal{A}\rightarrow\mathbb{R}_{\geq0}$ is a cost function that specifies the cost incurred by the team while navigating in the environment. The term $H(\pi(\cdot | s))=\sum_{a\in\mathcal{A}(s)}\pi(s,a)\log \pi(s,a)$ denotes the entropy of the policy $\pi$ in state $s\in\mathcal{S}$. Finally, the inefficiency parameter $\beta\in[0,\infty)$ balances the incurred costs with the randomness of the policy followed by the team.

The objective in \eqref{entropy_regularized_obj} corresponds to minimizing the entropy-regularized total cost where the weight of the regularization is controlled by the inefficiency parameter $\beta$. The constraints in \eqref{entropy_regularized_cons} ensure that the resulting behavior of the team in the environment satisfies the expected final distribution $\boldsymbol{\sigma}_{i}$. Note that as $\beta\rightarrow0$, the team is expected to reach its final distribution only through optimal paths that minimize their total cost. On the other hand, as $\beta\rightarrow\infty$, the team is expected to be as random as possible while reaching their final distribution. 

Let $\mathcal{G}\cup\mathcal{S}_0\cup\mathcal{S}_r$ be a partition of the set $\mathcal{S}$ where $\mathcal{G}$ is the set of goal states, $\mathcal{S}_0$ is the set of states from which there is no path reaching the states in $\mathcal{G}$, and $\mathcal{S}_r=\mathcal{S}\backslash \{\mathcal{G}\cup\mathcal{S}_0\}$. These sets can be efficiently computed via simple graph search algorithms, e.g., breadth-first search. By slightly modifying the results presented in \cite{savas2019entropy}, it can be shown that the problem in \crefrange{entropy_regularized_obj}{entropy_regularized_cons} is equivalent to the following convex optimization problem:
\small
\begin{subequations}
\begin{flalign}\label{entropy_regularized_start}
    &\underset{\substack{x(s,a)\geq 0}}{\text{minimize}} \; \sum_{s\in \mathcal{S}_r}\sum_{a\in \mathcal{A}(s)}x(s,a)\left[c(s,a) + \beta \log\left(\frac{x(s,a)}{\nu(s)}\right)\right] \\ 
    &\text{subject to:} \nonumber \\  \label{entropy_regularized_cons_1}
    & \nu(s)-\sum_{s'\in \mathcal{S}}\eta(s',s)= \alpha(s), \  \text{for all} \ s\in \mathcal{S}_r,\\  \label{entropy_regularized_cons_2}
    &\sum_{s\in \mathcal{S}_r}\eta(s,g_i) = \boldsymbol{\sigma}_{i}(g_i), \ \text{for all} \  g_i \in \mathcal{G},\\ \label{entropy_regularized_cons_3}
    &\eta(s,s') = \sum_{a\in \mathcal{A}(s)}x(s,a)P_{s,a,s'}, \ \text{for all} \  s \in \mathcal{S}_r  \ \text{and} \ s' \in \mathcal{S}, \\ 
    &\nu(s) = \sum_{a\in \mathcal{A}(s)}x(s,a), \ \text{for all} \ s \in \mathcal{S}_r.\label{entropy_regularized_end}
\end{flalign}
\end{subequations}
\normalsize
In the above problem, the decision variables $x(s,a)$ represent the density of the team members that occupy the state $s$ and take the action $a$. These variables are traditionally referred to as occupancy measures \cite{Altman}. The constraint in \eqref{entropy_regularized_cons_1} corresponds to balance equations, which express that the density entering a state should be equal to the density leaving that state. Similarly, the constraint in \eqref{entropy_regularized_cons_2} ensures that the final distribution of the team satisfies the condition in \eqref{entropy_regularized_cons}. Finally, the constraints in \crefrange{entropy_regularized_cons_3}{entropy_regularized_end} are introduced just to simplify the notation.

The objective in \eqref{entropy_regularized_start} is a convex function of $x(s,a)$ which combines linear terms $x(s,a)c(s,a)$ with the relative entropy of the distribution $x(s,a)$ with $\nu(s)$ \cite{boyd2004convex}. Since the constraints are also linear functions of $x(s,a)$, the resulting convex optimization problem can be solved efficiently via off-the-shelf solvers. However, to ensure the existence of optimal solutions, we need to choose the cost function $c(s,a)$ in a particular way, as described in the following proposition. 

{\setlength{\parindent}{0cm}
\noindent \begin{prop}
If the problem in \crefrange{opt_main_1}{opt_main_2} is feasible and $c(s,a)\geq\beta \log(\lvert \mathcal{A}(s)\rvert)$ for all $s\in\mathcal{S}_{r}$ and $a\in\mathcal{A}(s)$, then the the problem in \crefrange{entropy_regularized_start}{entropy_regularized_end} has a finite optimal solution. \label{suff_prop_for_finiteness}
\end{prop}}
\noindent\textbf{Proof:} We first note that if the problem in \eqref{opt_main_1}-\eqref{opt_main_2} is feasible, then the problem in \eqref{entropy_regularized_start}-\eqref{entropy_regularized_end} also has a feasible solution as shown in Lemma 1 in \cite{Marta}. Additionally, this feasible solution has a finite value due to the conventions that $0\log0$$=$$0$ and $0\log(0/0)$$=$$0$ which are based on continuity arguments. 

We now show that the optimal value is lower bounded by a finite constant. Suppose that $c(s,a)$$\geq$$\beta \log(\lvert \mathcal{A}(s)\rvert)$ for all $s$$\in$$\mathcal{S}$ and $a$$\in$$\mathcal{A}(s)$. For notational convenience, we drop the dependence of the set $\mathcal{A}(s)$ on $s$ in the following derivations. We express the objective function in \eqref{entropy_regularized_start} as $\sum_{s\in\mathcal{S}_r}\theta(s)$ where
\begin{align*}
    \theta(s):=\sum_{a\in\mathcal{A}}x(s,a)c(s,a) + \beta \sum_{a\in\mathcal{A}}x(s,a)\log\Big(\frac{x(s,a)}{\nu(s)}\Big).
\end{align*}
If $\nu(s)$$=$$0$, then we have $\theta(s)$$=$$0$ by the convention $0\log(0/0)$$=$$0$. Additionally, for each $s$$\in$$\mathcal{S}_r$ that satisfies $\nu(s)$$>$$0$, we have 
\begin{subequations}
\begin{flalign}\label{nonnegative_derivation_1}
    \theta(s)&\geq \beta\Big[\sum_{a\in\mathcal{A}}x(s,a)\log(\lvert \mathcal{A}\rvert) + \sum_{a\in\mathcal{A}}x(s,a)\log\Big(\frac{x(s,a)}{\nu(s)}\Big)\Big]&& \raisetag{20pt}\\ \label{nonnegative_derivation_2}
    & = \beta\Big[\nu(s)\log(\lvert \mathcal{A}\rvert) + \nu(s)\sum_{a\in\mathcal{A}}\frac{x(s,a)}{\nu(s)}\log\Big(\frac{x(s,a)}{\nu(s)}\Big)\Big]&&\\ \label{nonnegative_derivation_3}
    &\geq \beta \nu(s) \Big[\log(\lvert \mathcal{A}\rvert) - \log(\lvert \mathcal{A}\rvert)\Big] \geq 0.&& 
\end{flalign}
\end{subequations}
The inequality in \eqref{nonnegative_derivation_1} follows from the fact that $c(s,a)$$\geq$$\beta \log(\lvert \mathcal{A}(s)\rvert)$. The equality in \eqref{nonnegative_derivation_2} follows from the definition of $\nu(s)$ in \eqref{entropy_regularized_end} and the fact that $\nu(s)$$>$$0$. Finally, the inequality in \eqref{nonnegative_derivation_3} follows from the fact that the maximum entropy of a discrete probability distribution with a support size $K$$\in$$\mathbb{N}$ is always less than or equal to $\log(K)$.

Finally, since the problem in \eqref{entropy_regularized_start}-\eqref{entropy_regularized_end} has a feasible solution and its optimal value is lower bounded by zero, we conclude that it has a finite optimal solution. $\Box$

We note that if \(\beta\) is large, i.e., the agent's behavior is highly inefficient, then the \(x(s,a)\) may be unbounded. In this case, the agents are expected to spend infinite time in the environment thereby making the goal-directedness ineffective. The condition  $c(s,a)\geq\beta \log(\lvert \mathcal{A}(s)\rvert)$ given Proposition~\ref{suff_prop_for_finiteness} ensures that this pathological case does not happen, and inefficiency and goal-directedness are balanced.

Using the condition given in Proposition \ref{suff_prop_for_finiteness}, we can choose a cost function $c$ with sufficiently high values so that the problem in \crefrange{entropy_regularized_start}{entropy_regularized_end} has a finite solution. Let $\left\{x^{\star}(s,a): s\in\mathcal{S}, a\in\mathcal{A}(s)\right\}$ be a set of optimal decision variables for the problem in \crefrange{entropy_regularized_start}{entropy_regularized_end}. We can obtain the policy $\overline{\pi}$, which describes the expected behavior for the team to achieve the final distribution $\boldsymbol{\sigma}_{i}$, by the rule

\begin{align}\label{policy_construction}
    \overline{\pi}_i(s,a) = \begin{cases}
    \frac{x^{\star}(s,a)}{\sum_{a\in\mathcal{A}(s)}x^{\star}(s,a)} & \text{if} \ \ \sum_{a\in\mathcal{A}(s)}x^{\star}(s,a)>0\\
    \frac{1}{\lvert \mathcal{A}(s)\rvert} & \text{otherwise}.
    \end{cases}
\end{align}

In the following section, we will show how to utilize the policies $\overline{\pi}_i$ for quantifying deception and generating behaviors that manipulate the predictions of the adversary.

\section{Generating Deceptive Behavior} \label{sec:dec}
In this section, we introduce several measures to quantify the deceptiveness of the team's behavior and present efficient algorithms to synthesize deceptive policies. 

\subsection{Quantifying Deceptiveness through Statistical Distance}
We propose to quantify deception through the statistical distance between the team's behavior and the behavior expected by the adversary. Specifically, we utilize the Kullback-Leibler (KL) divergence to formally define deception.

{\setlength{\parindent}{0cm}
\noindent \begin{definition}
Let $Q_1$ and $Q_2$ be discrete probability distributions with a countable support $\mathcal{X}$. The Kullback-Leibler (KL) divergence between $Q_1$ and $Q_2$ is defined as 
\begin{align*}
    \text{KL}(Q_1 || Q_2) := \sum_{x\in\mathcal{X}} Q_1(x) \log \left(\frac{Q_1(x)}{Q_2(x)}\right).
\end{align*}
\end{definition}}

The KL divergence $\text{KL}(Q_1 || Q_2)$ measures the deviation of the distribution $Q_1$ from the distribution $Q_2$. As we will discuss shortly, in the context of deception, it provides us with a method to quantify the statistical deviation of the team's observed behavior from the behavior expected by the adversary. 

We consider two different types of deception, namely, exaggeration and ambiguity. Before providing formal definitions for these deceptive behaviors, we first introduce some notation. For an arbitrary policy $\pi\in\Pi(\mathbb{M})$, let $\Gamma^{\pi}$ be the distribution of paths in $\mathbb{M}$ generated under $\pi$. Note that the support of the distribution $\Gamma^{\pi}$ is the set $Paths_{\mathbb{M}}^{\pi,\alpha}$ of all paths, which may, in general, contain infinitely many elements. As we will shortly observe, for the purposes of deception, it is not necessary to explicitly construct this distribution.  

\textbf{Exaggeration:} In this first type of deception, the team aims to exaggerate its behavior to convince the adversary that they allocate their resources with respect to a \textit{decoy} goal distribution $\boldsymbol{\sigma}_{i}\in\mathcal{U}\setminus \{\boldsymbol{\sigma}_{*}\}$.
Without loss of generality, let $\boldsymbol{\sigma}_{1}$ be the true goal distribution, i.e., $\boldsymbol{\sigma}_{*}=\boldsymbol{\sigma}_{1}$. Then, for a given policy $\pi\in\Pi(\mathbb{M})$, we quantify the exaggeration of the team's resulting behavior through the following formula
\begin{align}\label{exagg_KL}
    &\textit{Deceptiveness}(\mathcal{U},\pi) = \max_{i\in [N]} \Big[ \text{KL}(\Gamma^{\pi}|| \Gamma^{\overline{\pi}_1}) - \text{KL}(\Gamma^{\pi}|| \Gamma^{\overline{\pi}_i})\Big].
\end{align}

The term $ \text{KL}(\Gamma^{\pi}|| \Gamma^{\overline{\pi}_i})$ quantifies the KL-divergence between the path distributions induced by the policies $\pi$ and $\overline{\pi}_i$. Therefore, the above deceptiveness metric measures the relative statistical distance of the paths induced by $\pi$ to the true policy $\overline{\pi}_1$ and decoy policy $\overline{\pi}_i$.

The intuition behind \eqref{exagg_KL} comes from the likelihood-ratio test, which is the most powerful hypothesis testing method for a given significance level \cite{neyman1933ix}. Recall that, for each $i\in[N]$, the adversary expects the team to follow a policy $\overline{\pi}_i$ to achieve the final distribution $\boldsymbol{\sigma}_{i}$. Now, suppose that the adversary runs the likelihood-ratio test to decide whether the team follows the policy $\overline{\pi}_i$ or $\overline{\pi}_j$. Let $\varrho_1,\ldots,\varrho_n$ be the paths followed by $n$ members of the team under the policy $\pi$. Moreover, let $\text{Pr}(\varrho_1,\ldots,\varrho_n | \overline{\pi}_i)$ and $\text{Pr}(\varrho_1,\ldots,\varrho_n | \overline{\pi}_j)$ be the probabilities of $\varrho_1,\ldots,\varrho_n$ under $\overline{\pi}_i$ and $\overline{\pi}_j$, respectively. By the likelihood-ratio test, for a given constant $C\in\mathbb{R}_{\geq 0}$, the adversary decides that the team aims to achieve the final distribution $\boldsymbol{\sigma}_{U^i}^{\star}$ through the policy $\overline{\pi}_i$ if 
\begin{align*}
       \log\Big(\text{Pr}(\varrho_1,\ldots,\varrho_n | \overline{\pi}_i)\Big)-\log\Big(\text{Pr}(\varrho_1,\ldots,\varrho_n | \overline{\pi}_j)\Big) \geq C.
\end{align*}

To see how \eqref{exagg_KL} is related to the likelihood-ratio test, note that, \(n\left[ \text{KL}(\Gamma^{\pi}|| \Gamma^{\overline{\pi}_1}) - \text{KL}(\Gamma^{\pi}|| \Gamma^{\overline{\pi}_i})\right]\) is equal to

\small
\begin{equation*}
      \mathbb{E}^{\pi}\Big[\log\Big(\text{Pr}(\varrho_1,\ldots,\varrho_n | \overline{\pi}_i)\Big)\Big]-\mathbb{E}^{\pi}\Big[\log\Big(\text{Pr}(\varrho_1,\ldots,\varrho_n | \overline{\pi}_1)\Big)\Big].
\end{equation*}
\normalsize

Therefore, the term inside the parenthesis in \eqref{exagg_KL} quantifies the expected log-likelihood of a goal distribution $\boldsymbol{\sigma}_{i}$ being the true goal distribution to the goal distribution $\boldsymbol{\sigma}_{1}$ being the true goal distribution when the team follows the policy $\pi$. Note that by taking the maximum over $i\in[N]$ in $\textit{Deceptiveness}(\mathcal{U},\pi)$, we quantify deceptiveness with respect to the most likely decoy goal distribution. Consequently, the problem in \crefrange{opt_main_1}{opt_main_2} corresponds to synthesizing a policy $\pi^{\star}$ that maximizes the expected relative log-likelihood for a decoy goal distribution $\boldsymbol{\sigma}_{i}\in\mathcal{U} \setminus\{\boldsymbol{\sigma}_{*}\}$ to be the true goal distribution while guaranteeing that the team's resulting behavior satisfies the final resource distribution $\boldsymbol{\sigma}_{*}$.

\textbf{Ambiguity:} In this second type of deception, the team aims to behave in a way to make its true goal distribution $\boldsymbol{\sigma}_{*}$ ambiguous to the adversary. Specifically, for a given policy $\pi\in\Pi(\mathbb{M})$, we quantify the ambiguity of the team's behavior through the following formula
\begin{align}\label{ambiguity_def}
    &\textit{Deceptiveness}(\mathcal{U},\pi) =  -\max_{i\in [N]}  \text{KL}(\Gamma^{\pi}|| \Gamma^{\overline{\pi}_i}).&&
\end{align}
Similar to the exaggeration behavior, the intuition behind the equation in \eqref{ambiguity_def} comes from the likelihood-ratio test. Specifically, in \eqref{ambiguity_def}, we measure the deceptiveness of a policy as the minimum expected log-likelihood of any utility matrix $\boldsymbol{\sigma}_{i}\in\mathcal{U}$. As a result, the problem in \crefrange{opt_main_1}{opt_main_2} corresponds to synthesizing a policy $\pi^{\star}$ that \textit{minimizes} the maximum log-likelihood for any goal distribution to be the true goal distribution while guaranteeing that the team's resulting behavior satisfies the final resource distribution $\boldsymbol{\sigma}_{*}\in\mathcal{U}$.

\subsection{Synthesis of Policies through Convex Optimization}
In this section, we present algorithms to solve the problem in \crefrange{opt_main_1}{opt_main_2} when $\textit{Deceptiveness}(\mathcal{U},\pi)$ is defined as in \eqref{exagg_KL} and in \eqref{ambiguity_def}. 

Although the problem in \crefrange{opt_main_1}{opt_main_2} is feasible, it is, in general, possible that the optimal value is not bounded below when $\textit{Deceptiveness}(\mathcal{U},\pi)$ is defined as in \eqref{exagg_KL} and in \eqref{ambiguity_def}. This is due to the fact that, for given $U^i$ and $U^j$, the support of the final distributions $\boldsymbol{\sigma}^{\star}_{1,U^i}$ and $\boldsymbol{\sigma}^{\star}_{1,U^j}$ may be different. As a result, the KL divergence between the path distributions $\Gamma^{\overline{\pi}_1}$ and $\Gamma^{\overline{\pi}_i}$ may be infinite. %

To ensure the finiteness of the optimal value in \crefrange{opt_main_1}{opt_main_2}, we propose to divide the team's behavior into two phases, namely, \textit{deceptive} and \textit{goal-directed phases}. During the deceptive phase, the team aims to deceive the adversary regarding its goal distribution by optimizing its behavior with respect to the measures in \eqref{exagg_KL} or in \eqref{ambiguity_def}. Let $T\in\mathbb{N}$ be a critical decision stage at which the team switches from the deceptive phase to the goal-directed phase. After $T$, the team aims to reach its final distribution $\boldsymbol{\sigma}_{*}$ through the shortest path.

We utilize extended MDPs to compactly represent the deceptive and goal-directed phases in a single decision model. Formally, let $\overline{\mathbb{M}}_T=(\overline{\mathcal{S}},\overline{\alpha},\mathcal{A}, \overline{P})$ denote an \textit{extended MDP} where $\overline{\mathcal{S}}=\mathcal{S}\times [T+1]$ is a finite set of states, $\overline{\alpha}:\overline{S}\rightarrow[0,1]$ is an initial distribution such that, for each $\langle s,t \rangle\in\overline{\mathcal{S}}$, $\overline{\alpha}(\langle s,t \rangle)=\alpha(s)$ if $t=1$ and $\overline{\alpha}(\langle s,t\rangle)=0$ otherwise, and $\overline{P}:\overline{\mathcal{S}}\times\mathcal{A}\times \overline{\mathcal{S}}\rightarrow[0,1]$ is a transition function such that
\begin{align*}
    \overline{P}_{\langle s,t\rangle, a, \langle s',t'\rangle} = \begin{cases} P_{s,a,s'} & \text{if}\ t \leq T \ \text{and}\ t' = t+1\\
    P_{s,a,s'} & \text{if}\ t = T+1 \ \text{and}\ t' = t\\
    0 & \text{otherwise}.
    \end{cases}
\end{align*}
In an extended MDP, we can clearly distinguish the deceptive and goal-directed phases by defining the objectives separately for the states $\mathcal{S}\times[T]$ and $\mathcal{S}\times\{T+1\}$ as will be discussed shortly.

In the above construction, $T$ is a design variable that can be used to tune the duration of the team's deceptive behavior. One practical approach is to set $T$ as a function of the shortest path. Specifically, let $T_{\min}$ be the minimum expected time for the team to reach their final distribution $\boldsymbol{\sigma}^{\star}_{1,U^i}$. $T_{\min}$ can be computed by replacing the objective function in \eqref{entropy_regularized_start} with $\sum_{s\in\mathcal{S}_r}\sum_{a\in\mathcal{A}(s)}x(s,a)$. Then, we can simply set $T=k\lceil T_{\min}\rceil$ where $k\in\mathbb{N}$ determines the balance between suboptimality of the behavior and deception effort. 

\noindent\textbf{Exaggeration:} We achieve the exaggeration behavior for deception by solving $N$ separate linear programs (LPs). Let ${\bf{x}}^{\pi}$ be the vector of occupancy measures that correspond to the policy $\pi$ constructed through the formula in \eqref{policy_construction}. By simple algebraic manipulations, it can be shown \cite{karabag2021deception} that, for each $i\in[N]$, we have 

\small
\begin{flalign*}
  &\text{KL}(\Gamma^{\pi}|| \Gamma^{\overline{\pi}_1}) - \text{KL}(\Gamma^{\pi}|| \Gamma^{\overline{\pi}_i})= \sum_{s\in \mathcal{S}_r}\sum_{a\in\mathcal{A}} {\bf{x}}^{\pi}(s,a)\log\Bigg(\frac{\overline{\pi}_i(s,a)}{\overline{\pi}_1(s,a)}\Bigg).&&
\end{flalign*}
\normalsize
Note in the above equation that the logarithmic term is a constant and corresponds to a virtual reward that quantifies the statistical likelihood of the decoy goal distribution $\boldsymbol{\sigma}_{i}$ with respect to the true goal distribution $\boldsymbol{\sigma}_{1}$. The virtual reward may be infinite when there is a support mismatch between the policies $\overline{\pi}_1$  and $\overline{\pi}_i$. To ensure the finiteness of the virtual reward and avoid computational issues, we add a small constant $\epsilon$ to both the numerator and the denominator in the logarithmic term. Accordingly, for each $i\in[N]$, we consider the following LP:

\small
\begin{subequations}
\begin{flalign}\label{exaggaration_start}
    &\underset{\substack{x(\langle s, t\rangle,a)\geq 0}}{\text{maximize}} \sum_{\langle s,t\rangle\in \mathcal{S}\times[T]}\sum_{a\in \mathcal{A}}x(\langle s,t\rangle,a)\log\Bigg(\frac{\overline{\pi}_i(s,a)+\epsilon}{\overline{\pi}_1(s,a)+\epsilon}\Bigg) \nonumber\\
    &\qquad \qquad\quad  -\sum_{\langle s,t\rangle\in \mathcal{S}\times\{T+1\}}\sum_{a\in \mathcal{A}}x(\langle s,t\rangle,a) 
    \\ 
    &\text{subject to:} \nonumber \\  \label{exaggaration_cons_1}
    & \nu(\langle s, t\rangle)-\sum_{\langle s',t'\rangle \in \mathcal{S}\times [T]}\eta(\langle s',t'\rangle,
\langle s,t\rangle)= \overline{\alpha}(\langle s,t\rangle), \nonumber \\
&\qquad \qquad \qquad \qquad \qquad \qquad \text{for all} \ \langle s,t\rangle\in \mathcal{S}_r\times[T+1]\\  \label{exaggaration_cons_2}
    &\sum_{\langle s,t \rangle \in \mathcal{S}_r\times [T]}\sum_{t'\in [T]}\eta(\langle s,t \rangle ,\langle g_i, t'\rangle ) = \sigma_{1,U^i}^{\star}(g_i), \ \text{for all} \  g_i \in \mathcal{G}\\ \label{exaggaration_cons_3}
    &\eta(\langle s,t \rangle ,\langle s',t'\rangle) = \sum_{a\in \mathcal{A}}x(\langle s, t\rangle,a) \overline{P}_{\langle s,t \rangle ,a,\langle s',t'\rangle}, \nonumber\\
    & \qquad \qquad\text{for all} \  \langle s,t\rangle \in \mathcal{S}_r\times [T]  \ \text{and} \ \langle s',t'\rangle \in \mathcal{S}\times [T] \\ 
    &\nu(\langle s,t \rangle) = \sum_{a\in \mathcal{A}}x(\langle s,t \rangle ,a), \ \text{for all} \ \langle s,t\rangle \in \mathcal{S}_r\times [T]\label{exaggaration_end}
\end{flalign}
\end{subequations}
\normalsize

The objective function in the above LP consists of two terms that enable the team to perform a deception phase followed by a goal-directed phase. Specifically, the first sum in the objective corresponds to $\text{KL}(\Gamma^{\pi}|| \Gamma^{\overline{\pi}_1}) - \text{KL}(\Gamma^{\pi}|| \Gamma^{\overline{\pi}_i})$ on the extended state-space $\mathcal{S}\times[T]$. On the other hand, the second term ensures that after the deception phase, i.e., $T+1$, the team reaches its final distribution by minimizing their total residence time in the environment.

The constraints in \crefrange{exaggaration_cons_1}{exaggaration_end} are the same as the constraints in \crefrange{entropy_regularized_cons_1}{entropy_regularized_end} with a minor difference. Specifically, the constraints in \crefrange{exaggaration_cons_1}{exaggaration_end} are now defined over the extended MDP $\overline{\mathbb{M}}_T$ instead of the original MDP $\mathbb{M}$.

Now, for each $i\in[N]$, let $v^{\star}_i$ be the optimal value of the LP given in \crefrange{exaggaration_start}{exaggaration_end} and $i^{\star}\in\arg\max_{i\in [N]}v^{\star}_i$. Moreover, let $\{x^{\star}(\langle s,t \rangle ,a): \langle s,t \rangle\in\mathcal{S}\times [T+1], \ a\in\mathcal{A}\}$ be the set of optimal decision variables corresponding to the LP with the optimal value $v^{\star}_{i^{\star}}$. We obtain an optimal deceptive policy $\pi^{\star}\in\Pi(\overline{\mathbb{M}}_T)$ through the construction 
\begin{flalign}\label{opt_policy}
    &\pi^{\star}\left(\langle s,t \rangle,a\right) =\nonumber\\
    &\begin{cases}
        \frac{x^{\star}(\langle s,t \rangle ,a)}{\sum_{a\in\mathcal{A}}x^{\star}(\langle s,t \rangle,a)} & \text{if} \ \ \sum_{a\in\mathcal{A}(\langle s,t \rangle)}x^{\star}(\langle s,t\rangle,a)>0\\
        \frac{1}{\lvert \mathcal{A}(\langle s,t\rangle)\rvert} & \text{otherwise}.
    \end{cases}
\end{flalign}

It follows from the standard results in the MDP theory, e.g., \cite[Chapter 7]{Puterman}, the policy $\pi^{\star}$ ensures that the team reaches its desired final distribution $\boldsymbol{\sigma}_{1}$.

\noindent \textbf{Ambiguity:} We achieve ambiguous behavior for deception by solving a single convex optimization problem. Recall from \eqref{ambiguity_def} that the objective in this type of deception is to obtain a policy ${\pi}$ that has the minimum statistical distance to each potential policy $\overline{\pi}_i$. Accordingly, using the derivations from the exaggeration behavior, we consider the following convex program: 

\small
\begin{subequations}
\begin{flalign}\label{ambiguity_start}
    &\underset{\substack{z, x(\langle s, t\rangle,a)\geq 0}}{\text{minimize}} \sum_{\langle s,t\rangle\in \mathcal{S}\times\{T+1\}}\sum_{a\in \mathcal{A}}x(\langle s,t\rangle,a) + z
    \\ 
    &\text{subject to:} \ \ \ \ \text{\crefrange{exaggaration_cons_1}{exaggaration_end}}, \nonumber \\  \label{ambiguity_cons_1}
    & z \geq \sum_{\langle s,t\rangle\in \mathcal{S}\times[T]}\sum_{a\in \mathcal{A}}x(\langle s,t\rangle,a)\log\left(\frac{x(\langle s,t\rangle,a)/ \sum_{a} x(\langle s,t\rangle,a)}{\overline{\pi}_i(s,a)+\epsilon}\right)  \nonumber \\  
&\qquad \qquad \qquad \qquad \qquad \qquad \qquad \qquad \text{for all} \ i \in [N]
\end{flalign}
\end{subequations}
\normalsize

The objective in the above convex program consists of two terms similar to the exaggeration behavior. The first term corresponds to the goal-directed phase in which the team aims to reach its final distribution by minimizing their occupancy time in the environment. The second term, expressed by the scalar variable $z$, corresponds to the ambiguity of the behavior during the deception phase. In particular, through the constraint in \eqref{ambiguity_cons_1}, this variable quantifies the maximum KL-divergence of the path distribution induced by $\pi$ to the set of path distributions induced by $\overline{\pi}_i$ where $i\in[N]$. 

Finally, an optimal deceptive policy for achieving ambiguous behavior can be obtained from the optimal decision variables for the program in \crefrange{ambiguity_start}{ambiguity_cons_1} through the construction introduced in \eqref{opt_policy}.

\newcommand{\StaticObstacle}[2]{ \fill[] (#1+0,#2+0) rectangle  (#1+1,#2+1);}
\newcommand{\initialstate}[2]{ \fill[black!30!brown] (#1+0.1,#2+0.1) rectangle (#1+0.9,#2+0.9);}
\newcommand{\goalstate}[2]{ \fill[black!50!green] (#1+0.1,#2+0.1) rectangle (#1+0.9,#2+0.9);}
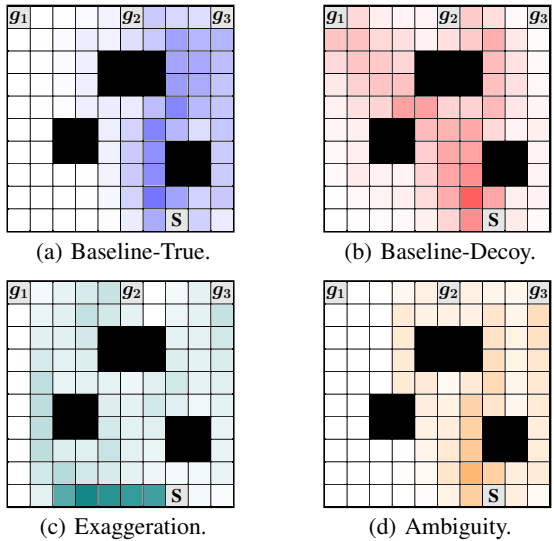
\begin{figure}[t!]
\centering
\begin{subfigure}[t]{0.23\textwidth}
\centering
\scalebox{0.3}{
\begin{tikzpicture}
\draw[black,line width=0.4pt] (0,0) grid[step=1, very thin] (10,10);
\draw[black,line width=3pt] (0,0) rectangle (10,10);

\fill[black!10!white] (10-0.05,0.05) rectangle (10-0.95,0.95);
\fill[black!10!white] (10-0.05,9.05) rectangle (10-0.95,9.95);
\fill[black!10!white] (10-4.05,9.05) rectangle (10-4.95,9.95);
\fill[black!10!white] (10-9.05,9.05) rectangle (10-9.95,9.95);
\node[] at (0.5,0.5) {\Huge \textbf{S}};
\node[] at (10-0.5,9.5) {\Huge $\boldsymbol{g_3}$};
\node[] at (10-4.5,9.5) {\Huge $\boldsymbol{g_2}$};
\node[] at (10-9.5,9.5) {\Huge $\boldsymbol{g_1}$};

\StaticObstacle{8}{2}
\StaticObstacle{7}{2}
\StaticObstacle{8}{3}
\StaticObstacle{7}{3}

\StaticObstacle{3}{3}
\StaticObstacle{2}{3}
\StaticObstacle{3}{4}
\StaticObstacle{2}{4}

\StaticObstacle{9-3}{6}
\StaticObstacle{9-4}{6}
\StaticObstacle{9-5}{6}
\StaticObstacle{9-3}{7}
\StaticObstacle{9-4}{7}
\StaticObstacle{9-5}{7}

\fill[blue!7.80!white] (10-0.05,0.05) rectangle (10-0.95,0.95);
\fill[blue!18.22!white] (10-1.05,0.05) rectangle (10-1.95,0.95);
\fill[blue!76.03!white] (10-2.05,0.05) rectangle (10-2.95,0.95);
\fill[blue!33.30!white] (10-3.05,0.05) rectangle (10-3.95,0.95);
\fill[blue!7.71!white] (10-4.05,0.05) rectangle (10-4.95,0.95);
\fill[blue!1.35!white] (10-5.05,0.05) rectangle (10-5.95,0.95);
\fill[blue!0.20!white] (10-6.05,0.05) rectangle (10-6.95,0.95);
\fill[blue!0.05!white] (10-7.05,0.05) rectangle (10-7.95,0.95);
\fill[blue!0.02!white] (10-8.05,0.05) rectangle (10-8.95,0.95);
\fill[blue!0.00!white] (10-9.05,0.05) rectangle (10-9.95,0.95);

\fill[blue!21.50!white] (10-0.05,1.05) rectangle (10-0.95,1.95);
\fill[blue!17.73!white] (10-1.05,1.05) rectangle (10-1.95,1.95);
\fill[blue!36.30!white] (10-2.05,1.05) rectangle (10-2.95,1.95);
\fill[blue!53.93!white] (10-3.05,1.05) rectangle (10-3.95,1.95);
\fill[blue!17.30!white] (10-4.05,1.05) rectangle (10-4.95,1.95);
\fill[blue!3.63!white] (10-5.05,1.05) rectangle (10-5.95,1.95);
\fill[blue!0.51!white] (10-6.05,1.05) rectangle (10-6.95,1.95);
\fill[blue!0.15!white] (10-7.05,1.05) rectangle (10-7.95,1.95);
\fill[blue!0.09!white] (10-8.05,1.05) rectangle (10-8.95,1.95);
\fill[blue!0.02!white] (10-9.05,1.05) rectangle (10-9.95,1.95);

\fill[blue!20.42!white] (10-0.05,2.05) rectangle (10-0.95,2.95);
\fill[blue!45.09!white] (10-3.05,2.05) rectangle (10-3.95,2.95);
\fill[blue!20.79!white] (10-4.05,2.05) rectangle (10-4.95,2.95);
\fill[blue!5.31!white] (10-5.05,2.05) rectangle (10-5.95,2.95);
\fill[blue!0.48!white] (10-6.05,2.05) rectangle (10-6.95,2.95);
\fill[blue!0.17!white] (10-7.05,2.05) rectangle (10-7.95,2.95);
\fill[blue!0.21!white] (10-8.05,2.05) rectangle (10-8.95,2.95);
\fill[blue!0.05!white] (10-9.05,2.05) rectangle (10-9.95,2.95);

\fill[blue!20.47!white] (10-0.05,3.05) rectangle (10-0.95,3.95);
\fill[blue!43.80!white] (10-3.05,3.05) rectangle (10-3.95,3.95);
\fill[blue!20.71!white] (10-4.05,3.05) rectangle (10-4.95,3.95);
\fill[blue!5.99!white] (10-5.05,3.05) rectangle (10-5.95,3.95);
\fill[blue!0.18!white] (10-8.05,3.05) rectangle (10-8.95,3.95);
\fill[blue!0.06!white] (10-9.05,3.05) rectangle (10-9.95,3.95);

\fill[blue!22.81!white] (10-0.05,4.05) rectangle (10-0.95,4.95);
\fill[blue!13.47!white] (10-1.05,4.05) rectangle (10-1.95,4.95);
\fill[blue!28.22!white] (10-2.05,4.05) rectangle (10-2.95,4.95);
\fill[blue!48.91!white] (10-3.05,4.05) rectangle (10-3.95,4.95);
\fill[blue!17.20!white] (10-4.05,4.05) rectangle (10-4.95,4.95);
\fill[blue!6.74!white] (10-5.05,4.05) rectangle (10-5.95,4.95);
\fill[blue!0.18!white] (10-8.05,4.05) rectangle (10-8.95,4.95);
\fill[blue!0.07!white] (10-9.05,4.05) rectangle (10-9.95,4.95);

\fill[blue!22.00!white] (10-0.05,5.05) rectangle (10-0.95,5.95);
\fill[blue!27.63!white] (10-1.05,5.05) rectangle (10-1.95,5.95);
\fill[blue!47.67!white] (10-2.05,5.05) rectangle (10-2.95,5.95);
\fill[blue!25.92!white] (10-3.05,5.05) rectangle (10-3.95,5.95);
\fill[blue!8.00!white] (10-4.05,5.05) rectangle (10-4.95,5.95);
\fill[blue!7.55!white] (10-5.05,5.05) rectangle (10-5.95,5.95);
\fill[blue!6.90!white] (10-6.05,5.05) rectangle (10-6.95,5.95);
\fill[blue!1.37!white] (10-7.05,5.05) rectangle (10-7.95,5.95);
\fill[blue!0.38!white] (10-8.05,5.05) rectangle (10-8.95,5.95);
\fill[blue!0.09!white] (10-9.05,5.05) rectangle (10-9.95,5.95);

\fill[blue!22.96!white] (10-0.05,6.05) rectangle (10-0.95,6.95);
\fill[blue!33.26!white] (10-1.05,6.05) rectangle (10-1.95,6.95);
\fill[blue!38.81!white] (10-2.05,6.05) rectangle (10-2.95,6.95);
\fill[blue!6.27!white] (10-6.05,6.05) rectangle (10-6.95,6.95);
\fill[blue!2.09!white] (10-7.05,6.05) rectangle (10-7.95,6.95);
\fill[blue!0.50!white] (10-8.05,6.05) rectangle (10-8.95,6.95);
\fill[blue!0.10!white] (10-9.05,6.05) rectangle (10-9.95,6.95);

\fill[blue!25.64!white] (10-0.05,7.05) rectangle (10-0.95,7.95);
\fill[blue!33.50!white] (10-1.05,7.05) rectangle (10-1.95,7.95);
\fill[blue!35.73!white] (10-2.05,7.05) rectangle (10-2.95,7.95);
\fill[blue!6.19!white] (10-6.05,7.05) rectangle (10-6.95,7.95);
\fill[blue!2.19!white] (10-7.05,7.05) rectangle (10-7.95,7.95);
\fill[blue!0.51!white] (10-8.05,7.05) rectangle (10-8.95,7.95);
\fill[blue!0.11!white] (10-9.05,7.05) rectangle (10-9.95,7.95);

\fill[blue!28.82!white] (10-0.05,8.05) rectangle (10-0.95,8.95);
\fill[blue!28.14!white] (10-1.05,8.05) rectangle (10-1.95,8.95);
\fill[blue!37.34!white] (10-2.05,8.05) rectangle (10-2.95,8.95);
\fill[blue!21.09!white] (10-3.05,8.05) rectangle (10-3.95,8.95);
\fill[blue!12.51!white] (10-4.05,8.05) rectangle (10-4.95,8.95);
\fill[blue!5.43!white] (10-5.05,8.05) rectangle (10-5.95,8.95);
\fill[blue!7.07!white] (10-6.05,8.05) rectangle (10-6.95,8.95);
\fill[blue!1.80!white] (10-7.05,8.05) rectangle (10-7.95,8.95);
\fill[blue!0.38!white] (10-8.05,8.05) rectangle (10-8.95,8.95);
\fill[blue!0.12!white] (10-9.05,8.05) rectangle (10-9.95,8.95);

\fill[blue!13.14!white] (10-1.05,9.05) rectangle (10-1.95,9.95);
\fill[blue!15.62!white] (10-2.05,9.05) rectangle (10-2.95,9.95);
\fill[blue!20.89!white] (10-3.05,9.05) rectangle (10-3.95,9.95);
\fill[blue!5.01!white] (10-5.05,9.05) rectangle (10-5.95,9.95);
\fill[blue!2.97!white] (10-6.05,9.05) rectangle (10-6.95,9.95);
\fill[blue!0.65!white] (10-7.05,9.05) rectangle (10-7.95,9.95);
\fill[blue!0.15!white] (10-8.05,9.05) rectangle (10-8.95,9.95);

\fill[black!10!white] (10-2.05,0.05) rectangle (10-2.95,0.95);
\node[] at (10-2.5,0.5) {\Huge \textbf{S}};

\end{tikzpicture}}
\vspace{-3pt} %
\caption{Baseline-True.} \label{truebeta6}
\vspace{5pt} %
\end{subfigure}
\begin{subfigure}[t]{0.23\textwidth}
\centering
\scalebox{0.3}{
\begin{tikzpicture}
\draw[black,line width=0.4pt] (0,0) grid[step=1, very thin] (10,10);
\draw[black,line width=3pt] (0,0) rectangle (10,10);

\fill[black!10!white] (10-0.05,0.05) rectangle (10-0.95,0.95);
\fill[black!10!white] (10-0.05,9.05) rectangle (10-0.95,9.95);
\fill[black!10!white] (10-4.05,9.05) rectangle (10-4.95,9.95);
\fill[black!10!white] (10-9.05,9.05) rectangle (10-9.95,9.95);
\node[] at (0.5,0.5) {\Huge \textbf{S}};
\node[] at (10-0.5,9.5) {\Huge $\boldsymbol{g_3}$};
\node[] at (10-4.5,9.5) {\Huge $\boldsymbol{g_2}$};
\node[] at (10-9.5,9.5) {\Huge $\boldsymbol{g_1}$};

\StaticObstacle{8}{2}
\StaticObstacle{7}{2}
\StaticObstacle{8}{3}
\StaticObstacle{7}{3}

\StaticObstacle{3}{3}
\StaticObstacle{2}{3}
\StaticObstacle{3}{4}
\StaticObstacle{2}{4}

\StaticObstacle{9-3}{6}
\StaticObstacle{9-4}{6}
\StaticObstacle{9-5}{6}
\StaticObstacle{9-3}{7}
\StaticObstacle{9-4}{7}
\StaticObstacle{9-5}{7}

\fill[red!2.19!white] (10-0.05,0.05) rectangle (10-0.95,0.95);
\fill[red!8.18!white] (10-1.05,0.05) rectangle (10-1.95,0.95);
\fill[red!76.03!white] (10-2.05,0.05) rectangle (10-2.95,0.95);
\fill[red!44.57!white] (10-3.05,0.05) rectangle (10-3.95,0.95);
\fill[red!16.50!white] (10-4.05,0.05) rectangle (10-4.95,0.95);
\fill[red!5.96!white] (10-5.05,0.05) rectangle (10-5.95,0.95);
\fill[red!2.63!white] (10-6.05,0.05) rectangle (10-6.95,0.95);
\fill[red!1.51!white] (10-7.05,0.05) rectangle (10-7.95,0.95);
\fill[red!0.74!white] (10-8.05,0.05) rectangle (10-8.95,0.95);
\fill[red!0.16!white] (10-9.05,0.05) rectangle (10-9.95,0.95);

\fill[red!5.67!white] (10-0.05,1.05) rectangle (10-0.95,1.95);
\fill[red!6.83!white] (10-1.05,1.05) rectangle (10-1.95,1.95);
\fill[red!35.09!white] (10-2.05,1.05) rectangle (10-2.95,1.95);
\fill[red!62.81!white] (10-3.05,1.05) rectangle (10-3.95,1.95);
\fill[red!32.64!white] (10-4.05,1.05) rectangle (10-4.95,1.95);
\fill[red!13.65!white] (10-5.05,1.05) rectangle (10-5.95,1.95);
\fill[red!6.26!white] (10-6.05,1.05) rectangle (10-6.95,1.95);
\fill[red!4.81!white] (10-7.05,1.05) rectangle (10-7.95,1.95);
\fill[red!3.41!white] (10-8.05,1.05) rectangle (10-8.95,1.95);
\fill[red!0.81!white] (10-9.05,1.05) rectangle (10-9.95,1.95);

\fill[red!5.25!white] (10-0.05,2.05) rectangle (10-0.95,2.95);
\fill[red!44.24!white] (10-3.05,2.05) rectangle (10-3.95,2.95);
\fill[red!35.31!white] (10-4.05,2.05) rectangle (10-4.95,2.95);
\fill[red!17.18!white] (10-5.05,2.05) rectangle (10-5.95,2.95);
\fill[red!5.17!white] (10-6.05,2.05) rectangle (10-6.95,2.95);
\fill[red!5.71!white] (10-7.05,2.05) rectangle (10-7.95,2.95);
\fill[red!8.06!white] (10-8.05,2.05) rectangle (10-8.95,2.95);
\fill[red!2.22!white] (10-9.05,2.05) rectangle (10-9.95,2.95);

\fill[red!5.25!white] (10-0.05,3.05) rectangle (10-0.95,3.95);
\fill[red!36.09!white] (10-3.05,3.05) rectangle (10-3.95,3.95);
\fill[red!33.90!white] (10-4.05,3.05) rectangle (10-4.95,3.95);
\fill[red!19.61!white] (10-5.05,3.05) rectangle (10-5.95,3.95);
\fill[red!6.75!white] (10-8.05,3.05) rectangle (10-8.95,3.95);
\fill[red!3.13!white] (10-9.05,3.05) rectangle (10-9.95,3.95);

\fill[red!5.86!white] (10-0.05,4.05) rectangle (10-0.95,4.95);
\fill[red!5.15!white] (10-1.05,4.05) rectangle (10-1.95,4.95);
\fill[red!15.17!white] (10-2.05,4.05) rectangle (10-2.95,4.95);
\fill[red!34.10!white] (10-3.05,4.05) rectangle (10-3.95,4.95);
\fill[red!29.43!white] (10-4.05,4.05) rectangle (10-4.95,4.95);
\fill[red!27.20!white] (10-5.05,4.05) rectangle (10-5.95,4.95);
\fill[red!6.44!white] (10-8.05,4.05) rectangle (10-8.95,4.95);
\fill[red!3.88!white] (10-9.05,4.05) rectangle (10-9.95,4.95);

\fill[red!4.75!white] (10-0.05,5.05) rectangle (10-0.95,5.95);
\fill[red!10.37!white] (10-1.05,5.05) rectangle (10-1.95,5.95);
\fill[red!27.35!white] (10-2.05,5.05) rectangle (10-2.95,5.95);
\fill[red!17.69!white] (10-3.05,5.05) rectangle (10-3.95,5.95);
\fill[red!17.10!white] (10-4.05,5.05) rectangle (10-4.95,5.95);
\fill[red!37.05!white] (10-5.05,5.05) rectangle (10-5.95,5.95);
\fill[red!36.76!white] (10-6.05,5.05) rectangle (10-6.95,5.95);
\fill[red!16.82!white] (10-7.05,5.05) rectangle (10-7.95,5.95);
\fill[red!12.19!white] (10-8.05,5.05) rectangle (10-8.95,5.95);
\fill[red!5.71!white] (10-9.05,5.05) rectangle (10-9.95,5.95);

\fill[red!3.73!white] (10-0.05,6.05) rectangle (10-0.95,6.95);
\fill[red!12.22!white] (10-1.05,6.05) rectangle (10-1.95,6.95);
\fill[red!24.85!white] (10-2.05,6.05) rectangle (10-2.95,6.95);
\fill[red!24.44!white] (10-6.05,6.05) rectangle (10-6.95,6.95);
\fill[red!23.14!white] (10-7.05,6.05) rectangle (10-7.95,6.95);
\fill[red!18.42!white] (10-8.05,6.05) rectangle (10-8.95,6.95);
\fill[red!9.19!white] (10-9.05,6.05) rectangle (10-9.95,6.95);

\fill[red!2.65!white] (10-0.05,7.05) rectangle (10-0.95,7.95);
\fill[red!11.92!white] (10-1.05,7.05) rectangle (10-1.95,7.95);
\fill[red!26.16!white] (10-2.05,7.05) rectangle (10-2.95,7.95);
\fill[red!16.09!white] (10-6.05,7.05) rectangle (10-6.95,7.95);
\fill[red!21.28!white] (10-7.05,7.05) rectangle (10-7.95,7.95);
\fill[red!23.10!white] (10-8.05,7.05) rectangle (10-8.95,7.95);
\fill[red!15.00!white] (10-9.05,7.05) rectangle (10-9.95,7.95);

\fill[red!1.33!white] (10-0.05,8.05) rectangle (10-0.95,8.95);
\fill[red!9.23!white] (10-1.05,8.05) rectangle (10-1.95,8.95);
\fill[red!31.21!white] (10-2.05,8.05) rectangle (10-2.95,8.95);
\fill[red!20.79!white] (10-3.05,8.05) rectangle (10-3.95,8.95);
\fill[red!12.55!white] (10-4.05,8.05) rectangle (10-4.95,8.95);
\fill[red!5.74!white] (10-5.05,8.05) rectangle (10-5.95,8.95);
\fill[red!12.02!white] (10-6.05,8.05) rectangle (10-6.95,8.95);
\fill[red!15.34!white] (10-7.05,8.05) rectangle (10-7.95,8.95);
\fill[red!23.87!white] (10-8.05,8.05) rectangle (10-8.95,8.95);
\fill[red!22.80!white] (10-9.05,8.05) rectangle (10-9.95,8.95);

\fill[red!2.81!white] (10-1.05,9.05) rectangle (10-1.95,9.95);
\fill[red!13.17!white] (10-2.05,9.05) rectangle (10-2.95,9.95);
\fill[red!20.70!white] (10-3.05,9.05) rectangle (10-3.95,9.95);
\fill[red!5.14!white] (10-5.05,9.05) rectangle (10-5.95,9.95);
\fill[red!4.49!white] (10-6.05,9.05) rectangle (10-6.95,9.95);
\fill[red!6.52!white] (10-7.05,9.05) rectangle (10-7.95,9.95);
\fill[red!15.30!white] (10-8.05,9.05) rectangle (10-8.95,9.95);

\fill[black!10!white] (10-2.05,0.05) rectangle (10-2.95,0.95);
\node[] at (10-2.5,0.5) {\Huge \textbf{S}};

\end{tikzpicture}}
\vspace{-3pt} %
\caption{Baseline-Decoy.} \label{decoybeta6}
\vspace{-5pt} %
\end{subfigure}

\begin{subfigure}[t]{0.23\textwidth}
\centering
\scalebox{0.3}{
\begin{tikzpicture}
\draw[black,line width=0.4pt] (0,0) grid[step=1, very thin] (10,10);
\draw[black,line width=3pt] (0,0) rectangle (10,10);

\fill[black!10!white] (10-0.05,0.05) rectangle (10-0.95,0.95);
\fill[black!10!white] (10-0.05,9.05) rectangle (10-0.95,9.95);
\fill[black!10!white] (10-4.05,9.05) rectangle (10-4.95,9.95);
\fill[black!10!white] (10-9.05,9.05) rectangle (10-9.95,9.95);
\node[] at (0.5,0.5) {\Huge \textbf{S}};
\node[] at (10-0.5,9.5) {\Huge $\boldsymbol{g_3}$};
\node[] at (10-4.5,9.5) {\Huge $\boldsymbol{g_2}$};
\node[] at (10-9.5,9.5) {\Huge $\boldsymbol{g_1}$};

\StaticObstacle{8}{2}
\StaticObstacle{7}{2}
\StaticObstacle{8}{3}
\StaticObstacle{7}{3}

\StaticObstacle{3}{3}
\StaticObstacle{2}{3}
\StaticObstacle{3}{4}
\StaticObstacle{2}{4}

\StaticObstacle{9-3}{6}
\StaticObstacle{9-4}{6}
\StaticObstacle{9-5}{6}
\StaticObstacle{9-3}{7}
\StaticObstacle{9-4}{7}
\StaticObstacle{9-5}{7}

\fill[teal!2.67!white] (10-0.05,0.05) rectangle (10-0.95,0.95);
\fill[teal!4.27!white] (10-1.05,0.05) rectangle (10-1.95,0.95);
\fill[teal!72.47!white] (10-2.05,0.05) rectangle (10-2.95,0.95);
\fill[teal!75.00!white] (10-3.05,0.05) rectangle (10-3.95,0.95);
\fill[teal!78.67!white] (10-4.05,0.05) rectangle (10-4.95,0.95);
\fill[teal!83.87!white] (10-5.05,0.05) rectangle (10-5.95,0.95);
\fill[teal!92.87!white] (10-6.05,0.05) rectangle (10-6.95,0.95);
\fill[teal!66.67!white] (10-7.05,0.05) rectangle (10-7.95,0.95);
\fill[teal!11.00!white] (10-8.05,0.05) rectangle (10-8.95,0.95);
\fill[teal!0.07!white] (10-9.05,0.05) rectangle (10-9.95,0.95);

\fill[teal!8.47!white] (10-0.05,1.05) rectangle (10-0.95,1.95);
\fill[teal!5.80!white] (10-1.05,1.05) rectangle (10-1.95,1.95);
\fill[teal!4.27!white] (10-2.05,1.05) rectangle (10-2.95,1.95);
\fill[teal!7.13!white] (10-3.05,1.05) rectangle (10-3.95,1.95);
\fill[teal!10.53!white] (10-4.05,1.05) rectangle (10-4.95,1.95);
\fill[teal!15.93!white] (10-5.05,1.05) rectangle (10-5.95,1.95);
\fill[teal!18.13!white] (10-6.05,1.05) rectangle (10-6.95,1.95);
\fill[teal!29.47!white] (10-7.05,1.05) rectangle (10-7.95,1.95);
\fill[teal!17.73!white] (10-8.05,1.05) rectangle (10-8.95,1.95);
\fill[teal!0.07!white] (10-9.05,1.05) rectangle (10-9.95,1.95);

\fill[teal!8.47!white] (10-0.05,2.05) rectangle (10-0.95,2.95);
\fill[teal!8.13!white] (10-3.05,2.05) rectangle (10-3.95,2.95);
\fill[teal!11.13!white] (10-4.05,2.05) rectangle (10-4.95,2.95);
\fill[teal!22.07!white] (10-5.05,2.05) rectangle (10-5.95,2.95);
\fill[teal!13.07!white] (10-6.05,2.05) rectangle (10-6.95,2.95);
\fill[teal!13.60!white] (10-7.05,2.05) rectangle (10-7.95,2.95);
\fill[teal!25.60!white] (10-8.05,2.05) rectangle (10-8.95,2.95);
\fill[teal!0.07!white] (10-9.05,2.05) rectangle (10-9.95,2.95);

\fill[teal!8.47!white] (10-0.05,3.05) rectangle (10-0.95,3.95);
\fill[teal!12.20!white] (10-3.05,3.05) rectangle (10-3.95,3.95);
\fill[teal!11.27!white] (10-4.05,3.05) rectangle (10-4.95,3.95);
\fill[teal!16.93!white] (10-5.05,3.05) rectangle (10-5.95,3.95);
\fill[teal!25.60!white] (10-8.05,3.05) rectangle (10-8.95,3.95);
\fill[teal!0.07!white] (10-9.05,3.05) rectangle (10-9.95,3.95);

\fill[teal!10.20!white] (10-0.05,4.05) rectangle (10-0.95,4.95);
\fill[teal!6.27!white] (10-1.05,4.05) rectangle (10-1.95,4.95);
\fill[teal!12.20!white] (10-2.05,4.05) rectangle (10-2.95,4.95);
\fill[teal!18.00!white] (10-3.05,4.05) rectangle (10-3.95,4.95);
\fill[teal!10.40!white] (10-4.05,4.05) rectangle (10-4.95,4.95);
\fill[teal!13.13!white] (10-5.05,4.05) rectangle (10-5.95,4.95);
\fill[teal!25.60!white] (10-8.05,4.05) rectangle (10-8.95,4.95);
\fill[teal!0.07!white] (10-9.05,4.05) rectangle (10-9.95,4.95);

\fill[teal!12.53!white] (10-0.05,5.05) rectangle (10-0.95,5.95);
\fill[teal!9.47!white] (10-1.05,5.05) rectangle (10-1.95,5.95);
\fill[teal!18.47!white] (10-2.05,5.05) rectangle (10-2.95,5.95);
\fill[teal!12.53!white] (10-3.05,5.05) rectangle (10-3.95,5.95);
\fill[teal!6.73!white] (10-4.05,5.05) rectangle (10-4.95,5.95);
\fill[teal!9.93!white] (10-5.05,5.05) rectangle (10-5.95,5.95);
\fill[teal!10.20!white] (10-6.05,5.05) rectangle (10-6.95,5.95);
\fill[teal!8.87!white] (10-7.05,5.05) rectangle (10-7.95,5.95);
\fill[teal!25.53!white] (10-8.05,5.05) rectangle (10-8.95,5.95);
\fill[teal!0.13!white] (10-9.05,5.05) rectangle (10-9.95,5.95);

\fill[teal!15.13!white] (10-0.05,6.05) rectangle (10-0.95,6.95);
\fill[teal!10.47!white] (10-1.05,6.05) rectangle (10-1.95,6.95);
\fill[teal!13.60!white] (10-2.05,6.05) rectangle (10-2.95,6.95);
\fill[teal!13.33!white] (10-6.05,6.05) rectangle (10-6.95,6.95);
\fill[teal!11.07!white] (10-7.05,6.05) rectangle (10-7.95,6.95);
\fill[teal!16.67!white] (10-8.05,6.05) rectangle (10-8.95,6.95);
\fill[teal!0.13!white] (10-9.05,6.05) rectangle (10-9.95,6.95);

\fill[teal!18.33!white] (10-0.05,7.05) rectangle (10-0.95,7.95);
\fill[teal!10.87!white] (10-1.05,7.05) rectangle (10-1.95,7.95);
\fill[teal!10.27!white] (10-2.05,7.05) rectangle (10-2.95,7.95);
\fill[teal!17.40!white] (10-6.05,7.05) rectangle (10-6.95,7.95);
\fill[teal!11.80!white] (10-7.05,7.05) rectangle (10-7.95,7.95);
\fill[teal!12.00!white] (10-8.05,7.05) rectangle (10-8.95,7.95);
\fill[teal!0.13!white] (10-9.05,7.05) rectangle (10-9.95,7.95);

\fill[teal!22.87!white] (10-0.05,8.05) rectangle (10-0.95,8.95);
\fill[teal!10.73!white] (10-1.05,8.05) rectangle (10-1.95,8.95);
\fill[teal!7.27!white] (10-2.05,8.05) rectangle (10-2.95,8.95);
\fill[teal!0.00!white] (10-3.05,8.05) rectangle (10-3.95,8.95);
\fill[teal!10.73!white] (10-4.05,8.05) rectangle (10-4.95,8.95);
\fill[teal!17.27!white] (10-5.05,8.05) rectangle (10-5.95,8.95);
\fill[teal!23.93!white] (10-6.05,8.05) rectangle (10-6.95,8.95);
\fill[teal!11.53!white] (10-7.05,8.05) rectangle (10-7.95,8.95);
\fill[teal!8.13!white] (10-8.05,8.05) rectangle (10-8.95,8.95);
\fill[teal!0.13!white] (10-9.05,8.05) rectangle (10-9.95,8.95);

\fill[teal!10.40!white] (10-1.05,9.05) rectangle (10-1.95,9.95);
\fill[teal!4.20!white] (10-2.05,9.05) rectangle (10-2.95,9.95);
\fill[teal!0.00!white] (10-3.05,9.05) rectangle (10-3.95,9.95);
\fill[teal!22.53!white] (10-5.05,9.05) rectangle (10-5.95,9.95);
\fill[teal!16.00!white] (10-6.05,9.05) rectangle (10-6.95,9.95);
\fill[teal!9.33!white] (10-7.05,9.05) rectangle (10-7.95,9.95);
\fill[teal!4.33!white] (10-8.05,9.05) rectangle (10-8.95,9.95);

\fill[black!10!white] (10-2.05,0.05) rectangle (10-2.95,0.95);
\node[] at (10-2.5,0.5) {\Huge \textbf{S}};

\end{tikzpicture}}
\vspace{-3pt} %
\caption{Exaggeration.} \label{exaggbeta6}
\end{subfigure}
\begin{subfigure}[t]{0.23\textwidth}
\centering
\scalebox{0.3}{
\begin{tikzpicture}
\draw[black,line width=0.4pt] (0,0) grid[step=1, very thin] (10,10);
\draw[black,line width=3pt] (0,0) rectangle (10,10);

\fill[black!10!white] (10-0.05,0.05) rectangle (10-0.95,0.95);
\fill[black!10!white] (10-0.05,9.05) rectangle (10-0.95,9.95);
\fill[black!10!white] (10-4.05,9.05) rectangle (10-4.95,9.95);
\fill[black!10!white] (10-9.05,9.05) rectangle (10-9.95,9.95);
\node[] at (0.5,0.5) {\Huge \textbf{S}};
\node[] at (10-0.5,9.5) {\Huge $\boldsymbol{g_3}$};
\node[] at (10-4.5,9.5) {\Huge $\boldsymbol{g_2}$};
\node[] at (10-9.5,9.5) {\Huge $\boldsymbol{g_1}$};

\StaticObstacle{8}{2}
\StaticObstacle{7}{2}
\StaticObstacle{8}{3}
\StaticObstacle{7}{3}

\StaticObstacle{3}{3}
\StaticObstacle{2}{3}
\StaticObstacle{3}{4}
\StaticObstacle{2}{4}

\StaticObstacle{9-3}{6}
\StaticObstacle{9-4}{6}
\StaticObstacle{9-5}{6}
\StaticObstacle{9-3}{7}
\StaticObstacle{9-4}{7}
\StaticObstacle{9-5}{7}

\fill[orange!4.80!white] (10-0.05,0.05) rectangle (10-0.95,0.95);
\fill[orange!12.40!white] (10-1.05,0.05) rectangle (10-1.95,0.95);
\fill[orange!75.87!white] (10-2.05,0.05) rectangle (10-2.95,0.95);
\fill[orange!37.87!white] (10-3.05,0.05) rectangle (10-3.95,0.95);
\fill[orange!9.80!white] (10-4.05,0.05) rectangle (10-4.95,0.95);
\fill[orange!2.47!white] (10-5.05,0.05) rectangle (10-5.95,0.95);
\fill[orange!0.40!white] (10-6.05,0.05) rectangle (10-6.95,0.95);
\fill[orange!0.13!white] (10-7.05,0.05) rectangle (10-7.95,0.95);
\fill[orange!0.07!white] (10-8.05,0.05) rectangle (10-8.95,0.95);
\fill[orange!0.00!white] (10-9.05,0.05) rectangle (10-9.95,0.95);

\fill[orange!15.27!white] (10-0.05,1.05) rectangle (10-0.95,1.95);
\fill[orange!13.67!white] (10-1.05,1.05) rectangle (10-1.95,1.95);
\fill[orange!36.20!white] (10-2.05,1.05) rectangle (10-2.95,1.95);
\fill[orange!55.73!white] (10-3.05,1.05) rectangle (10-3.95,1.95);
\fill[orange!20.20!white] (10-4.05,1.05) rectangle (10-4.95,1.95);
\fill[orange!6.07!white] (10-5.05,1.05) rectangle (10-5.95,1.95);
\fill[orange!0.93!white] (10-6.05,1.05) rectangle (10-6.95,1.95);
\fill[orange!0.07!white] (10-7.05,1.05) rectangle (10-7.95,1.95);
\fill[orange!0.07!white] (10-8.05,1.05) rectangle (10-8.95,1.95);
\fill[orange!0.07!white] (10-9.05,1.05) rectangle (10-9.95,1.95);

\fill[orange!14.60!white] (10-0.05,2.05) rectangle (10-0.95,2.95);
\fill[orange!41.40!white] (10-3.05,2.05) rectangle (10-3.95,2.95);
\fill[orange!17.73!white] (10-4.05,2.05) rectangle (10-4.95,2.95);
\fill[orange!8.80!white] (10-5.05,2.05) rectangle (10-5.95,2.95);
\fill[orange!0.60!white] (10-6.05,2.05) rectangle (10-6.95,2.95);
\fill[orange!0.07!white] (10-7.05,2.05) rectangle (10-7.95,2.95);
\fill[orange!0.07!white] (10-8.05,2.05) rectangle (10-8.95,2.95);
\fill[orange!0.07!white] (10-9.05,2.05) rectangle (10-9.95,2.95);

\fill[orange!14.33!white] (10-0.05,3.05) rectangle (10-0.95,3.95);
\fill[orange!36.00!white] (10-3.05,3.05) rectangle (10-3.95,3.95);
\fill[orange!13.60!white] (10-4.05,3.05) rectangle (10-4.95,3.95);
\fill[orange!11.00!white] (10-5.05,3.05) rectangle (10-5.95,3.95);
\fill[orange!0.07!white] (10-8.05,3.05) rectangle (10-8.95,3.95);
\fill[orange!0.07!white] (10-9.05,3.05) rectangle (10-9.95,3.95);

\fill[orange!15.87!white] (10-0.05,4.05) rectangle (10-0.95,4.95);
\fill[orange!6.67!white] (10-1.05,4.05) rectangle (10-1.95,4.95);
\fill[orange!19.27!white] (10-2.05,4.05) rectangle (10-2.95,4.95);
\fill[orange!33.47!white] (10-3.05,4.05) rectangle (10-3.95,4.95);
\fill[orange!9.67!white] (10-4.05,4.05) rectangle (10-4.95,4.95);
\fill[orange!13.27!white] (10-5.05,4.05) rectangle (10-5.95,4.95);
\fill[orange!0.07!white] (10-8.05,4.05) rectangle (10-8.95,4.95);
\fill[orange!0.07!white] (10-9.05,4.05) rectangle (10-9.95,4.95);

\fill[orange!17.67!white] (10-0.05,5.05) rectangle (10-0.95,5.95);
\fill[orange!8.93!white] (10-1.05,5.05) rectangle (10-1.95,5.95);
\fill[orange!29.20!white] (10-2.05,5.05) rectangle (10-2.95,5.95);
\fill[orange!16.60!white] (10-3.05,5.05) rectangle (10-3.95,5.95);
\fill[orange!5.47!white] (10-4.05,5.05) rectangle (10-4.95,5.95);
\fill[orange!16.33!white] (10-5.05,5.05) rectangle (10-5.95,5.95);
\fill[orange!16.33!white] (10-6.05,5.05) rectangle (10-6.95,5.95);
\fill[orange!0.00!white] (10-7.05,5.05) rectangle (10-7.95,5.95);
\fill[orange!0.07!white] (10-8.05,5.05) rectangle (10-8.95,5.95);
\fill[orange!0.07!white] (10-9.05,5.05) rectangle (10-9.95,5.95);

\fill[orange!19.87!white] (10-0.05,6.05) rectangle (10-0.95,6.95);
\fill[orange!9.87!white] (10-1.05,6.05) rectangle (10-1.95,6.95);
\fill[orange!25.40!white] (10-2.05,6.05) rectangle (10-2.95,6.95);
\fill[orange!16.33!white] (10-6.05,6.05) rectangle (10-6.95,6.95);
\fill[orange!0.00!white] (10-7.05,6.05) rectangle (10-7.95,6.95);
\fill[orange!0.07!white] (10-8.05,6.05) rectangle (10-8.95,6.95);
\fill[orange!0.07!white] (10-9.05,6.05) rectangle (10-9.95,6.95);

\fill[orange!22.53!white] (10-0.05,7.05) rectangle (10-0.95,7.95);
\fill[orange!9.87!white] (10-1.05,7.05) rectangle (10-1.95,7.95);
\fill[orange!22.67!white] (10-2.05,7.05) rectangle (10-2.95,7.95);
\fill[orange!16.33!white] (10-6.05,7.05) rectangle (10-6.95,7.95);
\fill[orange!0.00!white] (10-7.05,7.05) rectangle (10-7.95,7.95);
\fill[orange!0.07!white] (10-8.05,7.05) rectangle (10-8.95,7.95);
\fill[orange!0.07!white] (10-9.05,7.05) rectangle (10-9.95,7.95);

\fill[orange!26.07!white] (10-0.05,8.05) rectangle (10-0.95,8.95);
\fill[orange!9.07!white] (10-1.05,8.05) rectangle (10-1.95,8.95);
\fill[orange!20.47!white] (10-2.05,8.05) rectangle (10-2.95,8.95);
\fill[orange!10.67!white] (10-3.05,8.05) rectangle (10-3.95,8.95);
\fill[orange!10.47!white] (10-4.05,8.05) rectangle (10-4.95,8.95);
\fill[orange!9.80!white] (10-5.05,8.05) rectangle (10-5.95,8.95);
\fill[orange!16.33!white] (10-6.05,8.05) rectangle (10-6.95,8.95);
\fill[orange!0.00!white] (10-7.05,8.05) rectangle (10-7.95,8.95);
\fill[orange!0.07!white] (10-8.05,8.05) rectangle (10-8.95,8.95);
\fill[orange!0.13!white] (10-9.05,8.05) rectangle (10-9.95,8.95);

\fill[orange!7.20!white] (10-1.05,9.05) rectangle (10-1.95,9.95);
\fill[orange!7.93!white] (10-2.05,9.05) rectangle (10-2.95,9.95);
\fill[orange!11.47!white] (10-3.05,9.05) rectangle (10-3.95,9.95);
\fill[orange!11.33!white] (10-5.05,9.05) rectangle (10-5.95,9.95);
\fill[orange!6.47!white] (10-6.05,9.05) rectangle (10-6.95,9.95);
\fill[orange!0.00!white] (10-7.05,9.05) rectangle (10-7.95,9.95);
\fill[orange!0.00!white] (10-8.05,9.05) rectangle (10-8.95,9.95);

\fill[black!10!white] (10-2.05,0.05) rectangle (10-2.95,0.95);
\node[] at (10-2.5,0.5) {\Huge \textbf{S}};

\end{tikzpicture}}
\vspace{-3pt} %
\caption{Ambiguity.} \label{ambibeta6}
\end{subfigure}

\caption{Normalized density distributions of the autonomous team for \(\beta = 6\) and $c(s,a)$$=$$10$. The true final distribution is $[0.5, 0.5, 0]$ and the decoy final distribution is $[0, 0.5, 0.5]$ (a) The baseline density distribution (blue) expected by the adversary for the true final distribution. (b) The baseline density distribution (red) expected by the adversary for the decoy final distribution. (c) The density distribution (teal) for the exaggeration behavior. (d) The density distribution (orange) for the ambiguity behavior.}
\vspace{-10pt} %
\label{tunable_figure}
\end{figure}

\section{Simulations}
In this section, we present a numerical simulation to illustrate the proposed deception strategy in a motion planning example. We utilize the CVXPY interface \cite{diamond2016cvxpy} and the ECOS solver \cite{domahidi2013ecos} to obtain solutions to the considered convex optimization problems.  

We consider a scenario in which the team navigates in an environment represented by the $10\times10$ grid-world shown in \cref{tunable_figure}. Each grid cell represents a state, and there are four available actions $\{up,down,right, left\}$ in each state under which the agents transition to the neighboring state in the corresponding direction. All team members start their motion from the state labeled with $S$. The team aims to allocate their resources over three goal states labeled with $g_1$, $g_2$, and $g_3$. We consider two potential goal allocations $\boldsymbol{\sigma}_{1}=[0, 0.5, 0.5]$ and $\boldsymbol{\sigma}_{2}=[0.5, 0.5, 0]$. Note that since we have $\boldsymbol{\sigma}_{1}=[0,0.5,0.5]$, the team aims to allocate half of the team members to $g_2$ and the remaining half to $g_3$.  

We synthesize a deceptive density control strategy for Team 1 using the proposed methods. For the synthesis of deceptive control strategies, we solve the optimization problem \crefrange{exaggaration_start}{exaggaration_end} for exaggeration behavior and the optimization problem \crefrange{ambiguity_start}{ambiguity_cons_1} for ambiguity behavior. 

In \cref{tunable_figure}, we illustrate the density distributions that the adversary predicts, as well as the distribution that the team follows. The density of state \(s\) is equal to \(\sum_{t \in [T+1]} \sum_{a \in \mathcal{A}} x(\langle s,t\rangle, a)\) that is the expected time that the team members spend at state \(s\).

For the synthesis of deceptive strategies, we set the parameter $T=5$ in the optimization problems, i.e., the team shows the deceptive behavior for \(5\) steps and then switches to the goal-directed behavior. In \cref{exaggbeta6}, we observe that Team 1 exaggerates its behavior (shown in green color) and pretends to achieve the decoy distribution $[0.5,0.5,0]$ during the initial deception phase. Then, during the goal-directed phase, they eventually reach their final distribution. For the ambiguity behavior shown in \cref{ambibeta6}, during the initial deception phase, the team follows a path (shown in orange color) that is the only significantly plausible path for both the true and decoy final distributions. This behavior preserves the ambiguity of the true distribution until the goal-directed phase.

\section{User Study}
We now assess the effectiveness of the proposed deceptive strategies in altering real users' perceptions of final goal distribution by showcasing five paths sampled from each policy (see \cref{tunable_figure}) at steps $4, 8, 12, 16, 20$.

\subsection{Experiment Design}
{\bf Manipulated Variables:}
We conducted two sets of experiments. In the first set, participants were informed of the deceptive purpose (YES), while in the second set, the true purpose was not disclosed until after the study (NO).
Within each set, we varied the policy types: a baseline policy for the true goal distribution (\textit{Baseline-True}), a baseline policy for the decoy goal distribution (\textit{Baseline-Decoy}), a deceptive policy that creates ambiguity regarding the true goal distribution (\textit{Ambiguity}), or a deceptive policy that exaggerates towards a decoy goal distribution (\textit{Exaggeration}). The baseline policies are generated using the principle of maximum entropy method discussed in Section \ref{sec:maxentropy}, and the deceptive policies are generated using the methods given in \ref{sec:dec}. We use \(\beta = 6\) for the maximum entropy method.
In total, the experiment has $8$ different conditions.

{\bf Dependent Measures:}
We have two dependent measures: correctness, measured as the ratio of correctly perceived goals, and confidence on a 5-point Likert scale. We combine the two in a \textbf{score}: the confidence if they get two goals correct, half the confidence if they get only one goal correct, and negative of the confidence if they are not correct at all~\cite{huang2019enabling}. This score captures that if one is incorrect, it is better not to be confident about it. The lower the score is, the more effective a deceptive policy is.

{\bf Participants:}
We used a between-subjects design, where each participant would only see paths from one condition, in order to avoid carryover and fatigue effects from having seen a different condition before.
We ran this experiment on a total of $320$ participants across the $8$ conditions, recruited via Prolific \cite{palan2018prolific}. We filtered out data from $3$ ($0.938\%$) participants who withdrew consent.
The average age of the $317$ consenting participants was $39.595$ (SD = $13.939$). The gender ratio was $0.497$ female.

\subsection{Results}
Across every condition studied, we observed increases in both correctness and confidence as the path advanced. This outcome is as expected, as longer path segments generally reduce interpretive uncertainty for participants. By the final step, the path endpoints in all conditions either reach or nearly reach the predetermined goals, leading to high average and low variance in scores. \footnote{All data in the results is publicly available at \url{https://github.com/vivianchen98/deception_user_study_data}.}

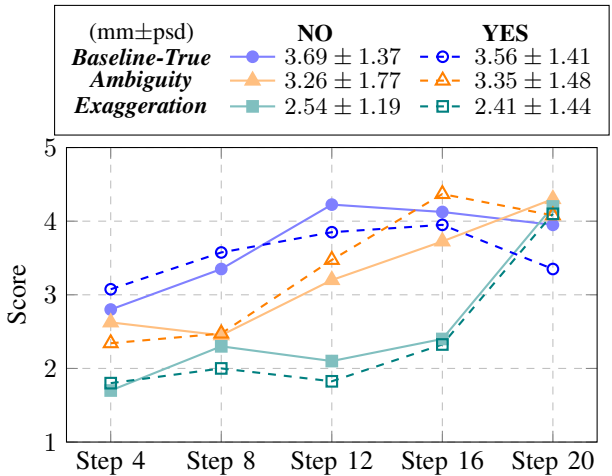
\begin{figure}
  \centering
      \begin{tikzpicture}
            \begin{axis}[
                width=\linewidth,
                height=5.5cm,
                ylabel={Score},
                grid=major,
                grid style=dashed,
                symbolic x coords={step4, step8, step12, step16, step20},
                xticklabels={0, Step 4, Step 8, Step 12, Step 16, Step 20},
                ymin=1, ymax=5,
                ytick={0,1,2,3,4,5},
            ]
            
                \addplot[
                    color=blue!50!white,
                    solid,
                    thick,
                    mark=*,
                ] table [x=step, y=prediction_0_scores_mean, y error=prediction_0_scores_std, col sep=comma] {exp_results/NO_6_scores.dat};
                \label{plot:NO_pred_true}
            
                \addplot[
                    color=orange!50!white,
                    solid,
                    thick,
                    mark=triangle*,
                    mark size=3pt,
                ] table [x=step, y=ambiguity_scores_mean, y error=ambiguity_scores_std, col sep=comma] {exp_results/NO_6_scores.dat};
                \label{plot:NO_ambi}
    
                \addplot[
                    color=teal!50!white,
                    solid,
                    thick,
                    mark=square*,
                ] table [x=step, y=exaggeration_scores_mean, y error=exaggeration_scores_std, col sep=comma] {exp_results/NO_6_scores.dat};
                \label{plot:NO_exagg}
    
                \addplot[
                    color=blue,
                    dashed,
                    thick,
                    mark=o,
                    mark options={
                        solid,
                        draw=blue, %
                        fill=white,      %
                    }
                ] table [x=step, y=prediction_0_scores_mean, y error=prediction_0_scores_std, col sep=comma] {exp_results/YES_6_scores.dat};
                \label{plot:YES_pred_true}
                
                \addplot[
                    color=orange,
                    dashed,
                    thick,
                    mark=triangle,
                    mark size=3pt,
                    mark options={
                        solid,
                        draw=orange, %
                    }
                ] table [x=step, y=ambiguity_scores_mean, y error=ambiguity_scores_std, col sep=comma] {exp_results/YES_6_scores.dat};
                \label{plot:YES_ambi}
            
                \addplot[
                    color=teal,
                    dashed,
                    thick,
                    mark=square,
                    mark options={
                        solid,
                        draw=teal, %
                    }
                ] table [x=step, y=exaggeration_scores_mean, y error=exaggeration_scores_std, col sep=comma] {exp_results/YES_6_scores.dat};
                \label{plot:YES_exagg}
                
            \end{axis}
    
            \matrix (m) [matrix of nodes, 
                         nodes in empty cells,
                         nodes={anchor=center},
                         draw,
                         font=\small,
                         column sep=0, row sep=-5,
                         anchor=north, %
                         at={(current axis.north)}, yshift=1.9cm
                        ] 
            {
                (mm$\pm$psd) & \textbf{NO} & \textbf{YES}\\
                \textbf{\textit{Baseline-True}} & \ref{plot:NO_pred_true} $3.69\pm 1.37$ & \ref{plot:YES_pred_true} $3.56\pm 1.41$\\
                \textbf{\textit{Ambiguity}} & \ref{plot:NO_ambi} $3.26\pm 1.77$ & \ref{plot:YES_ambi} $3.35\pm 1.48$\\
                \textbf{\textit{Exaggeration}} & \ref{plot:NO_exagg} $2.54\pm 1.19$ & \ref{plot:YES_exagg} $2.41\pm 1.44$\\
            };
            
        \end{tikzpicture}
        \caption{Participant scores for paths sampled from the \textit{Baseline-True}, \textit{Ambiguity}, and \textit{Exaggeration} policies, in both NO and YES groups (all with $\beta=6$). The legend delineates statistical details for each line, including the mean of means (mm) and the pooled standard deviation (psd) across path steps.}
        \label{fig:three_plots}
        \vspace{-10pt} %
\end{figure}

{\bf H1: The proposed deceptive policies are effective:}
\cref{fig:three_plots} plots the mean scores for the Baseline-True, Ambiguity, and Exaggeration policies at each step of the path within both the NO and YES groups. Whether in the YES or NO group, participants consistently register lower mean scores for the Ambiguity and Exaggeration policies compared to the Baseline-True policy, a pattern particularly pronounced before Step 16. Upon closer examination of the two deceptive policies, we notice that Exaggeration yields even lower mean scores than Ambiguity. This not only supports H1 but further reveals exaggeration is more effective than ambiguity.

{\bf H2: Informing participants of deception makes a policy less effective:}
We made this hypothesis based on the presumption that users are more likely to second-guess themselves if they know they are being deceived, leading to reduced scores. Nevertheless, the outcomes presented in \cref{fig:three_plots} from both the YES and NO groups for each policy type---depicted by solid and dashed lines of the same color---do not manifest a consistent trend, rendering the findings non-definitive with respect to this hypothesis.

\begin{figure}
 \centering
       \begin{tikzpicture}
            \begin{axis}[
                width=\linewidth,
                height=5.5cm,
                ylabel={Score},
                grid=major,
                grid style=dashed,
                symbolic x coords={step4, step8, step12, step16},
                xtick=data,
                xticklabels={Step 4, Step 8, Step 12, Step 16},
                ymin=1, ymax=5,
                ytick={0,1,2,3,4,5},
            ]
        
            \addplot[
                color=blue!50!white,
                solid, thick,
                mark=otimes*,
                mark options={solid, draw=blue!50!white, fill=white}
            ] table [x=step, y=prediction_0_scores_mean, y error=prediction_0_scores_std, col sep=comma] {exp_results/NO_6_scores_4step.dat};
            \label{plot:original_true}
        
            \addplot[
                color=red!50!white,
                solid, thick,
                mark=diamond*,
                mark size=3pt,
            ] table [x=step, y=prediction_1_scores_mean, y error=prediction_1_scores_std, col sep=comma] {exp_results/NO_6_scores_4step.dat};
            \label{plot:original_decoy}
    
            \addplot[
                color=blue,
                dashed, thick,
                mark=otimes,
                mark options={solid}
            ] table [x=step, y=prediction_0_scores_mean, y error=prediction_0_scores_std, col sep=comma] {exp_results/flipped_NO_6_scores.dat};
            \label{plot:flipped_true}
        
            \addplot[
                color=red,
                dashed, thick,
                mark=diamond,
                mark size=3pt,
                mark options={solid}
            ] table [x=step, y=prediction_1_scores_mean, y error=prediction_1_scores_std, col sep=comma] {exp_results/flipped_NO_6_scores.dat};
            \label{flipped_decoy}
            
            \end{axis}
    
        \matrix (m) [matrix of nodes, 
                     nodes in empty cells,
                     nodes={anchor=center},
                     draw,
                     font=\small,
                     column sep=0, row sep=-5,
                     anchor=north, %
                     at={(current axis.north)}, yshift=1.6cm
                    ] 
        {
             & \textbf{Original} & \textbf{Flipped}\\
            \textbf{\textit{Baseline-True}} ($g_2+g_3$) & \ref{plot:original_true} right & \ref{plot:flipped_true} left\\
            \textbf{\textit{Baseline-Decoy}} ($g_1+g_2$) & \ref{plot:original_decoy} left & \ref{flipped_decoy} right\\
        };
        \end{tikzpicture}

        \caption{Participant scores for paths sampled from the \textit{Baseline-True} and \textit{Baseline-Decoy} policies across path steps in the original and the flipped experiments (both with $\beta=6$ in NO group).}
        \vspace{-10pt} %
        \label{fig:baselines}
\end{figure}
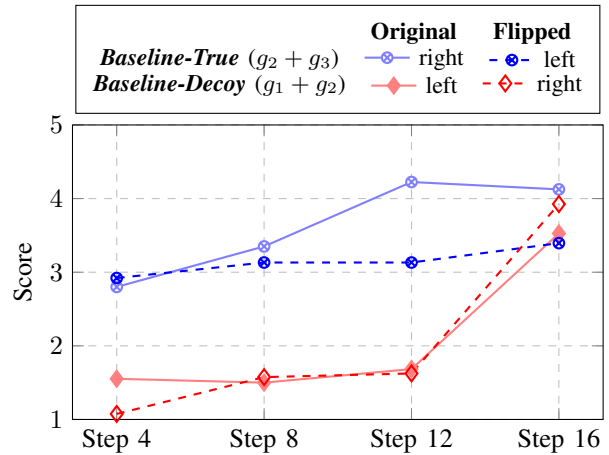

{\bf H3: Baseline-True and Baseline-Decoy have comparable scores:}
\cref{fig:baselines} shows that the Baseline-True policy, with goals closer to the starting position on the right side of the interface, consistently attains higher scores than Baseline-Decoy (see solid lines).
Two conceivable reasons may account for this trend: participants might either possess an inherent preference for goals located nearer to the starting point, or they could have a bias towards the right side of the interface, perhaps due to the ease of clicking.
To discern the actual bias, we mirrored the experiment using horizontally-flipped paths, using data from $315$ consenting participants. We limited the display to four steps, as the final step would not yield substantial insights. 
The results for the flipped experiment, represented by dashed lines in \cref{fig:baselines}, show a similar trend, suggesting the higher scores for Baseline-True are not a consequence of the goals' alignment on the interface. Rather, it reveals that participants have a consistent inclination towards proximate goals. 
This inclination can also be explained by the notion of last deceptive point (LDP) proposed \cite{masters2017deceptive}, as the Baseline-Decoy paths inadvertently approach the LDP whereas the Baseline-True paths never come close to the LDP.

\section{Conclusions}
We studied the problem of synthesizing deceptive resource allocation strategies for a team consisting of a large number of autonomous agents. We developed a prediction algorithm based on the principle of maximum entropy that models the predictions of adversarial observers regarding the autonomous team's final allocation strategy over multiple goal locations. By quantifying deceptiveness as a function of statistical distance between certain distributions, we then developed deceptive strategies, based on convex optimization, to control the density of the team members in the environment while they progress towards their final distribution. A user study validates the effectiveness of the proposed algorithms.

\addtolength{\textheight}{-12cm}   %

\bibliographystyle{IEEEtran}
\bibliography{main}

\end{document}